\documentclass[11pt,oneside,english,reqno]{amsart}
\usepackage[T1]{fontenc}
\usepackage[utf8]{inputenc}
\usepackage[a4paper]{geometry}
\usepackage{mathrsfs}
\usepackage{amstext}
\usepackage{amsthm}
\usepackage{amssymb}
\usepackage{bbm}
\usepackage{shuffle}
\usepackage{color}
\usepackage{float}
\usepackage{enumerate}
\usepackage{verbatim} 
\usepackage{stmaryrd}
\usepackage{enumitem}
\usepackage{graphicx}
\usepackage{todonotes}
\usepackage{cancel}
\usepackage[normalem]{ulem}
\usepackage{subcaption}

\usepackage{array}
\usepackage{tabularx}

\usepackage{hyperref}

\usetikzlibrary{arrows.meta, shapes.geometric}
\usepackage{tikz, pgfplots}
\usepackage{tikz-cd}
\usetikzlibrary{matrix,chains,positioning,decorations.pathreplacing,arrows}
\usetikzlibrary{calc,positioning}

\usepackage{graphicx}
\usepackage{caption}
\definecolor{diagramblue}{HTML}{315873}
\definecolor{diagramink}{HTML}{27333D}

\makeatletter
\numberwithin{equation}{section}
\numberwithin{figure}{section}
\theoremstyle{plain}
\newtheorem{thm}{\protect\theoremname}
\numberwithin{thm}{section}
\theoremstyle{definition}
\newtheorem{defn}[thm]{\protect\definitionname}
\theoremstyle{remark}
\newtheorem{rem}[thm]{\protect\remarkname}
\theoremstyle{plain}
\newtheorem{lem}[thm]{\protect\lemmaname}
\theoremstyle{plain}
\newtheorem{prop}[thm]{\protect\propositionname}
\theoremstyle{plain}

\theoremstyle{plain}
\newtheorem{example}[thm]{\protect\examplename}

\makeatother

\usepackage{babel}
\providecommand{\corollaryname}{Corollary}
\providecommand{\definitionname}{Definition}
\providecommand{\lemmaname}{Lemma}
\providecommand{\propositionname}{Proposition}
\providecommand{\remarkname}{Remark}
\providecommand{\theoremname}{Theorem}
\providecommand{\examplename}{Example}

\newcommand{\cA}{\mathcal{A}}

\newcommand{\cF}{\mathcal{F}}
\newcommand{\cG}{\mathcal{G}}
\newcommand{\cH}{\mathcal{H}}
\newcommand{\cI}{\mathcal{I}}

\newcommand{\cL}{\mathcal{L}}

\newcommand{\cP}{\mathcal{P}}

\newcommand{\cS}{\mathcal{S}}
\newcommand{\cT}{\mathcal{T}}

\newcommand{\cV}{\mathcal{V}}

\newcommand{\sig}{\mathbf{Sig}}
\newcommand{\bsig}{\mathbf{BSig}}
\newcommand{\vsig}{\mathbf{VSig}}
\newcommand{\gvsig}{\mathbf{GVSig}}

\newcommand{\geovsig}{\mathbf{gVSig}}
\newcommand{\trees}{\mathcal{T}}
\newcommand{\ktree}{\tau_{\kappa}}
\newcommand{\dtree}{\tau_{d}}
\newcommand{\vtree}{\tau_{\nu}}

\newcommand{\kforest}{f_{\kappa}}
\newcommand{\dforest}{f_{d}}
\newcommand{\vforest}{f_{\nu}}

\newcommand{\tset}{\mathcal{T}_{\nu}}
\newcommand{\fset}{\mathcal{F}_{\nu}}

\newcommand{\kprod}{\mathbin{{\star}_{\kappa}}}
\newcommand{\dprod}{\mathbin{{\star}_{d}}}
\newcommand{\vprod}{\mathbin{{\star}_{\nu}}}

\newcommand{\KK}{\mathbb{K}}

\newcommand{\NN}{\mathbb{N}}

\newcommand{\RR}{\mathbb{R}}

\newcommand{\VV}{\mathbb{V}}

\newcommand{\bA}{\mathbf{A}}

\newcommand{\bF}{\mathbf{F}}
\newcommand{\bG}{\mathbf{G}}

\newcommand{\bk}{\mathbf{k}}

\newcommand{\bR}{\mathbf{R}}

\newcommand{\bV}{\mathbf{V}}
\newcommand{\bW}{\mathbf{W}}
\newcommand{\bX}{\mathbf{X}}

\newcommand{\bY}{\mathbf{Y}}
\newcommand{\by}{\mathbf{y}}
\newcommand{\bZ}{\mathbf{Z}}

\newcommand{\se}{\mathsf e}

\newcommand{\bv}{\mathbf{v}}

\newcommand{\bw}{\mathbf{w}}

\newcommand{\hopf}{\mathscr{H}}

\newcommand{\kparent}[1]{p_{\tau_{\kappa}}(#1)}

\newcommand{\myspan}{\mbox{span}}

\newcommand{\Id}{\mathrm{Id}}

\newcommand{\dd}{\mathop{}\!\mathrm{d}}

\newenvironment{tree}
               {\begin{tikzpicture}[every node/.style={circle,draw,fill,minimum size=3pt,inner sep=0pt,outer sep=0pt},line cap=round,baseline = .1pt]}
               {\end{tikzpicture}}

\newcommand{\innerprod}[1]{\left\langle #1 \right\rangle}

\usepackage[external]{forest}

\def\dtR<#1>{\Forest{[#1]}}
\def\dtI<#1#2>{\Forest{[#1[#2]]}}
\def\dtII<#1#2#3>{\Forest{[#1[#2[#3]]]}}
\def\dtV<#1#2#3>{\Forest{[#1[#2][#3]]}}

\forestset{
  decor/.style={
    label/.expanded={
      [inner sep=0.1ex,font=\unexpanded{\tiny}]
      right:{$#1$}
    }
  },
  root/.style={
    minimum size=0.5ex
  },
  Decorated/.style={
    for tree={
      circle,
      fill,
      inner sep=0ex,
      minimum size=0.5ex,
      grow'=north,
      l=0,
      l sep=0.8ex,
      s sep=0.5em,
      fit=tight,
      parent anchor=center,
      child anchor=center,
      delay={
        decor/.option=content,
        content=
      }
    }
  },
  default preamble={Decorated,root},
  begin draw/.code={%
    \csname begin\endcsname{tikzpicture}[
      baseline={([yshift=-0.5ex]current bounding box.center)}
    ]%
  },
}

\title{The Generalized Volterra Signature}
\author{Munawar Ali$^*$\qquad Qi Feng$^\dagger$\qquad Samy Tindel$^\sharp$}

\thanks{$^*$: Department of
Mathematics, Florida State University, Tallahassee, 32306; email: ma22bm@fsu.edu.}
\thanks{$^\dagger$: Department of
Mathematics, Florida State University, Tallahassee, 32306; email: qfeng2@fsu.edu. This author is partially supported by the National Science Foundation under grant \#DMS-2420029.}
\thanks{$^\sharp$: Department of
Mathematics, Purdue University,  West Lafayette, 47907; email: stindel@purdue.edu.}
\date{\today}

\begin{document}
\sloppy
\maketitle

\begin{abstract}
In this paper, we introduce the generalized Volterra signature \(\gvsig(\bX)\) of a path \(\bX\), extending the classical Volterra signature introduced in [Hager, Harang, Pelizzari \& Tindel, arXiv:2603.04525, 2026], which extends Lyons' path signature to the Volterra setting. We establish a new algebraic structure to encode the generalized Volterra signature. We also establish a generalized Chen's identity. By restricting the differential trees to ladder trees, we introduce the geometric Volterra signature \(\geovsig(\bX)\) and establish the shuffle product of its signature components. To the best of our knowledge, the shuffle product for the Volterra signature is the first of its kind and allows us to formulate a universal approximation theorem for geometric Volterra signatures. For a special class of kernel functions, we prove the existence of an extended path \(\Bar{\bX}\) and relate its classical signature to the geometric Volterra signature of the original path \(\bX\).
\end{abstract}




\section{Introduction}
The dynamics of many real--world phenomena reflect an intricate dependence on history, with recent fluctuations and distant events jointly shaping their evolution. These phenomena are modeled via dynamical systems that encode memory via kernel function(s). There are myriad examples of such dynamical systems in literature. For example, these systems are used to model volatility of financial assets \cite{vl_ComteRenault1996,vl_Gatheral2018,vl_GuyonLekeufack2023}, long and short memory in engineering and signal--processing problems \cite{vl_BoydChua1985,vl_Leland1994,vl_BurovBarkai2008,vl_CarmonaCoulon2014}, viscoelastic materials in physics \cite{vl_DiethelmFreed1999}, spread of epidemics in biology \cite{vl_Cromer1985}, and biological and neural systems with memory \cite{vl_GerstnerKistler2002,vl_Kou2008}.

Machine learning offers several ways to study history--dependent behaviour from sequential data. One line of work builds trainable continuous--time models through neural ordinary differential equations \cite{vl_Chen2018} and neural controlled differential equations, the latter being designed to handle irregularly sampled time series \cite{vl_Kidger2020}. Attention-based architectures, including Transformers and their efficient variants, learn relationships between observations at different positions in a sequence \cite{vl_Vaswani2017,vl_Wang2020}. Recurrent neural networks (RNNs) summarize previous observations through successive updates of a hidden state \cite{vl_Elman1990,vl_Jordan1997}, while long short--term memory networks (LSTMs) use gates to control the retention and modification of this information \cite{vl_Hochreiter1997}. State-space architectures such as S4 and Mamba offer further mechanisms for processing long sequences through structured or selective state updates \cite{vl_Gu2022,vl_GuDao2023,vl_DaoGu2024}. The choice of architecture therefore determines how information from the past is stored, selected, and used to form predictions.

Despite their flexibility, these approaches present challenges in understanding and controlling the influence of past observations. The representations learned by neural models can be difficult to interpret, which may obscure the temporal dependencies underlying their predictions \cite{vl_Lipton2018,vl_Rudin2019}. In recurrent neural networks, learning over long time horizons can also be hindered by vanishing or exploding gradients \cite{vl_Bengio1994,vl_Pascanu2013}. For scientific applications, a further consideration is how to incorporate prior knowledge about the underlying dynamics directly into model design and training \cite{vl_Raissi2019,vl_Karniadakis2021}. These considerations motivate approaches that combine flexible modeling of temporal dependence with an explicit and mathematically tractable description of memory.

A classical approach to modeling such memory-dependent phenomena is through Volterra integral equations and their stochastic counterparts. A typical system of such equations takes the form
\begin{equation} \label{SVIE}
    \bY_t^i = y_i + \sum_{j=0}^{d} \int_0^t k_i(t,r)\,V_j^i(\bY_r)\,d \bX_r^j, \qquad i=1,\ldots,n,
\end{equation}
where \(\bY_t = (\bY_t^1,\ldots,\bY_t^n)\) is the response process, \(\bX^j_r\) are independent drivers for \(j=1,2,\dots,d\) with \(d \bX_r^0:=dr\), and \(k_i(t,r)\), for \(i=1,\ldots,n\) are kernel functions that encode the memory.

Such equations have a well-established theory and applications in control \cite{vl_Corduneanu2008,vl_Gripenberg1990} and stochastic analysis \cite{vl_PardouxProtter1990,vl_OksendalZhang1993}. They also provide a framework for learning dynamics from data, as illustrated by neural integral equations \cite{vl_Zappala2024} and neural stochastic Volterra equations \cite{vl_Bergerhausen2025}. These developments extend the use of neural differential equation models \cite{vl_Chen2018,vl_Kidger2020} to integral formulations that can explicitly incorporate memory.

A natural approach to studying \eqref{SVIE} is to identify the iterated integrals that govern its solution. For controlled differential equations this perspective is based on Chen's work on iterated integrals \cite{chen1977iterated} and Lyons' theory of rough paths \cite{lyons1998differential}, which provides a framework for interpreting differential equations driven by irregular signals. In the Volterra setting, the corresponding iterated integrals must also account for the interactions between integration times introduced by the kernels. Harang and Tindel \cite{harang2021volterra} developed a rough path framework for Volterra equations with singular kernels, subsequently extended by Harang, Tindel, and Wang to higher-order expansions \cite{harang2023volterra} and to probabilistic constructions of Volterra rough paths \cite{harang2024volterra}. Bruned and Katsetsiadis \cite{bruned2023ramification} developed a Hopf-algebraic framework for Volterra rough paths using decorated rooted trees, including planted trees, and drawing on ideas from branched rough path theory. More recently, Hager, Harang, Pelizzari, and Tindel \cite{hager2026volterra} developed \emph{the Volterra signature} as a feature representation for history-dependent systems. Their construction organizes the iterated integrals arising from linear controlled Volterra equations into a tensor-valued object that incorporates the memory kernel and satisfies a convolutional Chen identity.

Following the stochastic Taylor expansion approach of Feng and Zhang \cite{feng2023cubature}, we consider an expansion of the system of stochastic Volterra integral equations \eqref{SVIE}. To separate the upper integration limit from the time argument of the kernel, we introduce the two-time process
\begin{equation}\label{eq:two-time-process}
    \Phi_t^{i,\tau}
    :=
    y_i+\sum_{j=0}^{d}\int_0^t
    k_i(\tau,r)\,V_j^i(\bY_r)\,d\bX_r^j,
    \qquad 0\leq t\leq\tau\leq T,
\end{equation}
whose diagonal recovers the original solution:
$\Phi_t^{i,t}=\bY_t^i$. To approximate a smooth function $F$ of the terminal response $\bY_T$, we apply the Taylor expansion iteratively. The resulting terms involve derivatives of $F$ and the vector fields, together with iterated integrals weighted by products of kernels. Retaining products of iterated integrals in their formal product form, we obtain contributions such as the following:
\begin{align} \label{representative third order term}
    & F(\bY_T) \approx \, \mathrm{LOT} \nonumber \\
    & \quad + \sum_{i_1,i_2,i_3=1}^{n}\sum_{j_1,j_2,j_3=0}^{d}
    C_{i_1,i_2,i_3}^{j_1,j_2,j_3}(\by)
    \int_0^T k_{i_1}(T,r_1)
    \left(\int_0^{r_1}k_{i_2}(T,r_2)\,d\bX_{r_2}^{j_2}\right)
    \left(\int_0^{r_1}k_{i_3}(r_1,r_3)\,d\bX_{r_3}^{j_3}\right)
    d\bX_{r_1}^{j_1} \nonumber \\
    & + \mathrm{HOT},
\end{align}

where \(C_{i_1,i_2,i_3}^{j_1,j_2,j_3}(\by) := \partial_{i_2}\partial_{i_1}F(\by)\, \partial_{i_3}V_{j_1}^{i_1}(\by)\, V_{j_2}^{i_2}(\by)V_{j_3}^{i_3}(\by)\),  and \(\mathrm{LOT}\) denotes the lower-order terms, while \(\mathrm{HOT}\) collects all the remaining terms of degree three and higher.

The displayed term is naturally represented by the compatible Volterra tree


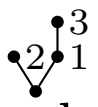
\[
\vtree = \left(
\Forest{Decorated[[2][1[3]]]},
\Forest{Decorated[1[2][3]]}
\right),
\]
where \(\Forest{Decorated[[2][1[3]]]}\) is the kernel tree denoted by \(\ktree\) and \(\Forest{Decorated[1[2][3]]}\) is the differential tree denoted by \(\dtree\). The dummy root in the kernel tree corresponds to the terminal time \(T\). The kernel tree \(\ktree\) encodes the product of the kernels
\[
    K_{i_1}(T,r_1)K_{i_2}(T,r_2)K_{i_3}(r_1,r_3),
\]
where the variables \(r_1\) and \(r_2\) are attached to the terminal variable \(T\), while \(r_3\) is attached to \(r_1\). On the other hand, the differential tree \(\dtree\) encodes the order of the differentials. In particular, \(\dtree = \Forest{Decorated[1[2][3]]}\) encodes the following object
\[
    \int_0^{T}\cdots\, \left(\int_0^{r_1}\cdots\,d\bX_{r_2}^{j_2}\right)
    \left(\int_0^{r_1}\cdots\,d\bX_{r_3}^{j_3}\right)
    d\bX_{r_1}^{j_1}.
\]
Thus, the kernel tree and the differential tree can have different structures. This provides a natural motivation for introducing \emph{compatible Volterra tree} as a pair consisting of a kernel tree and a differential tree, and for studying the algebraic structures associated with by such objects in the following sections.

Based on this new algebraic structure, we define the \emph{generalized Volterra signature}, denoted by \(\gvsig\), of a path \(\bX:[0,T]\to\RR^d\). A typical component of this signature takes the form
\begin{equation}
    \int_s^{t} k_{i_1}(\tau,r_1)
    \left(\int_0^{r_1}k_{i_2}(\tau,r_2)\,d\bX_{r_2}^{j_2}\right)
    \left(\int_0^{r_1}k_{i_3}(r_1,r_3)\,d\bX_{r_3}^{j_3}\right)
    d\bX_{r_1}^{j_1}.
\end{equation}
Here, \(\tau\) denotes the upper time parameter of the kernels, with \(0\leq s<t\leq\tau\leq T\). Setting \(s=0\) and \(t=\tau=T\) recovers the integral appearing in \eqref{representative third order term}.

Our main contributions can be summarized as follows:

\begin{itemize}
    \item \emph{Algebraic structure.} We introduce compatible Volterra trees, each consisting of a kernel tree and a differential tree that satisfy a compatibility condition. We define the compatible Volterra forests as the formal product of compatible Volterra trees and construct a Hopf algebra on these forests, encoding kernel interactions and the structure of iterated integrals separately. This framework provides the algebraic foundation for the generalized Volterra signature \(\gvsig\).
    \item \emph{The generalized convolution and Chen's identity.} We define a generalized convolution product for the components of generalized Volterra signature and establish the corresponding Chen's identity, describing how the generalized Volterra signature over an interval is recovered from its values on adjacent subintervals by incorporating all kernel interactions.
    \item \emph{Shuffle product and the weakly geometric Volterra signature \(\geovsig\).} We define a space of linear combinations of compatible Volterra trees whose differential trees are ladder trees and define a Volterra shuffle product on this space. We also prove that this space is closed under the Volterra shuffle product, i.e., the result is a linear combination of compatible Volterra trees whose differential trees remain ladder trees. The resulting shuffle identity along with Chen's identity and analytic bounds provide the basis for defining the weakly geometric Volterra signature \(\geovsig\).
    \item \emph{Injectivity and universal approximation.} We prove the injectivity of the generalized Volterra signature and the weakly geometric Volterra signature for time-extended \(\alpha\)-Hölder paths. Using the Stone--Weierstrass theorem, we establish universal approximation theorems for both signatures.
    \item \emph{Reduction of \(\geovsig\) to \(\sig\) of an extended path for a class of kernels.} For the special class of separable kernels, we prove that each component of the weakly geometric Volterra signature can be recovered from a suitable truncation of the classical signature of an extended path. 
\end{itemize}

\section{Preliminaries}
To make the manuscript self-contained, we introduce the notation and path spaces used throughout the paper. We then recall the definition and basic properties of the classical signature, introduce the classical Volterra signature, and review rooted trees and the Connes--Kreimer Hopf algebra of rooted trees.

\subsection{Notation and path spaces}
We begin by introducing some basic notation and specifying the path spaces considered in this paper. For the definitions and standard results recalled below, we refer to Friz and Victoir \cite{Friz_Victoir_2010}.

\begin{defn}[\(p\)-variation of a path]
Let \(p \geq 1\) and let \(\bX:[0,T]\to\RR^d\) be a path.
Its \(p\)-variation on \([0,T]\) is defined by
\[\|\bX\|_{p\text{-var};[0,T]} := \left( \sup_{\cP} \sum_{i=0}^{n-1} \big|\bX_{t_{i+1}}-\bX_{t_i}\big|^p \right)^{1/p},
\]
where the supremum is taken over all finite partitions
\(\cP=\{0=t_0<t_1<\cdots<t_n=T\}\).
\end{defn}
We also define the paths of bounded \(p\)-variation as follows.
\begin{defn}[Path of bounded \(p\)-variation]
A path \(\bX:[0,T] \to \RR^d\) is said to be a path of bounded \(p\)-variation if
\begin{equation*}
    \|\bX\|_{p\text{-var};[0,T]}  < \infty.
\end{equation*}
We denote the set of paths of bounded \(p\)-variation as \(\cV_p([0,T];\RR^d)\). Unless otherwise stated, we will be working with paths of bounded \(p\)-variation for \(p=1\), also called paths of bounded variation. 
\end{defn}
A related notion of path regularity is \(\alpha\)-H\"older continuity. We now define this notion and explain its relationship with finite \(p\)-variation.
\begin{defn}[\(\alpha\)-H\"older paths]
Let \(\alpha\in(0,1]\). A path \(\bX:[0,T]\to\RR^d\)
is said to be \(\alpha\)-H\"older continuous if
\begin{equation*}
    \|\bX\|_{\alpha\text{-H\"ol};[0,T]} := \sup_{0\le s<t\le T}
    \frac{|\bX_t-\bX_s|}{(t-s)^\alpha} <\infty.
\end{equation*}
The quantity \(\|\bX\|_{\alpha\text{-H\"ol};[0,T]}\)
is called the \(\alpha\)-H\"older seminorm of \(\bX\).
\end{defn}
The following lemma relates H\"older continuity to finite
\(p\)-variation.
\begin{lem}
Let \(\alpha\in(0,1]\). Every \(\alpha\)-H\"older continuous
path \(\bX:[0,T]\to\RR^d\) has finite \(p\)-variation
for every \(p\geq 1/\alpha\). Moreover,
\begin{equation*}
    \|\bX\|_{p\text{-var};[0,T]}
    \leq T^\alpha
    \|\bX\|_{\alpha\text{-H\"ol};[0,T]}.
\end{equation*}
\end{lem}

\subsection{Words and shuffle algebra}
For an \(\RR^d\)-valued path \(\bX\), we fix the alphabet set to be \(\cA = \{1,2,\dots, d\}\) and call this the set of letters. We define a \textbf{word} \(\bw\) of length \(n\) to be a sequence \(\bw = w_1 w_2 \dots w_n\) where \(w_i \in \cA\) for \(i = 1,2,\dots, n\). The length of the word \(\bw\) is the number of letters in it and is denoted by \(|\bw|\). Next, we define some operations on the words as follows.
\begin{defn}[Concatenation of words]
    The concatenation of two words \(\bv = v_1 \dots v_m\) and \(\bw=w_1\dots w_n\), denoted by \(\bv \bw\), is defined below.
    \[
    \bv \bw = v_1 \dots v_m w_1\dots w_n.
    \]
\end{defn}

\begin{defn}[Shuffle product of words]
    For two words \(\bv = v_1 \dots v_m\) and \(\bw=w_1\dots w_n\), we define their shuffle product \(\bv \shuffle \bw\) recursively as follows
    \[
    \bw \shuffle \emptyset  = \emptyset \shuffle \bw = \bw, \quad \text{ and } 
    \]
    \[
    \bw \shuffle \bv  = \left[(w_1\dots w_{n-1}) \shuffle \bv \right] w_n + \left[\bw \shuffle (v_1 \dots v_{m-1})\right]v_m,
    \]
    where \(\emptyset\) is an empty word.
\end{defn}
Moreover, we denote \(\bW\) to be the set of all words and \(\bW_n\) to be the set of all words of length \(n\). Also, for \(n=0\) we define \(\bW_0\) to be the set with empty word \(\emptyset\) . We denote the vector space generated by the set of words \(\bW\) as \(\VV (\cA)\) defined as follows
\[\VV (\cA) = \left\{ \ell = \sum_{\bw \in \bW} c_{\bw} \bw : c_{\bw} \in \RR , c_{\bw} \ne 0 \text{ for finitely many } \bw \in \bW\right\}.\]
\begin{defn}
The set \(\VV (\cA)\) endowed with product defined as the concatenation of words is a non-commutative algebra.
\end{defn}
We define the \textbf{extended tensor algebra} over \(\RR^d\) to be the space
\[T((\RR^d)) = \prod_{n=0}^{\infty}(\RR^d)^{\otimes n} = \left\{\bv=(v_0,v_1,\dots) \;\middle|\; v_n\in(\RR^d)^{\otimes n},\ n=0,1,\dots\right\},\]
where \((\RR^d)^{\otimes0}:=\RR\). Another space, closely related to this one and to the space generated by words, is the \textbf{tensor algebra} over \(\RR^d\), defined as
\[T(\RR^d) = \bigoplus_{n=0}^{\infty}(\RR^d)^{\otimes n} = \left\{\bv\in T((\RR^d)) \;\middle|\; \exists \; K\in\NN \text{ such that } v_n=0 \text{ for all } n\ge K\right\}.\]

Also, the \textbf{truncated tensor algebra} over \(\RR^d\) is defined as
\[T^N(\RR^d) := \left\{\bv\in T(\RR^d) \;\middle|\; v_n=0 \text{ for all } n>N\right\}.\]
Through the coordinate pairing, \(\VV(\cA)\) can be identified with a subspace of the algebraic dual of the extended tensor algebra \(T((\RR^d))\). That is,
\[\VV(\cA)\hookrightarrow T((\RR^d))^*.\]

\subsection{Signature and its properties}
On the interval \([0,T]\) we define the \(n\)-simplex as follows.
\[
\Delta^n_T : =\{(r_1,\dots,r_n): r_1 \le \dots \le r_n, \,  r_i \in [0,T], \forall i=1,\dots,n\}.
\]
We define the signature of a path of bounded variation \(\bX:[0,T] \to \RR^d\) as follows.
\begin{defn}
The signature of a continuous \(\RR^d\)-valued path of bounded variation \((\bX_t)_{t \in [0,T]}\) is the \(T((\RR^d))\)-valued process \( (s,t) \in \Delta_T^2 \mapsto \sig(\bX)_{st} \in T((\RR^d))\) whose component corresponding to each word \(\bw = w_1 \dots w_n \in \bW \subsetneq \VV (\cA), n \ge 1 \) is defined as
\[\innerprod{\sig(\bX)_{st}, \bw} := \int_{s}^{t} \int_{s}^{r_n} \cdots \int_{s}^{r_2} d\bX^{w_1}_{r_1}\cdots d\bX^{w_n}_{r_n},\]
where the integrals are defined in Riemann--Stieltjes sense. Similarly, for an empty word \(\emptyset\) we define \(\innerprod{\sig(\bX)_{st}, \emptyset} := 1\).
Equivalently, the signature can be written in terms of its tensor levels as
\[\sig(\bX)_{st} = \left(1, \int_s^t d\bX_r, \int_s^t \int_s^{r_2} d\bX_{r_1}\otimes d\bX_{r_2}, \int_s^t \int_s^{r_3}\int_s^{r_2} d\bX_{r_1}\otimes d\bX_{r_2}\otimes d\bX_{r_3}, \dots\right).\]
\end{defn}
Furthermore, the level \(N\) truncation of the signature is given as below
\[\sig^N(\bX)_{st} = \sum_{\bw \, \in \, \bW \,:\, |\bw|\, \le \, N} \innerprod{\sig(\bX)_{st}, \bw} e_{\bw},\]
where \(e_{\bw}\) is a basis element of \(T((\RR^d))\).
\begin{rem}
Throughout, for 
\(s=0\) we employ the following convention \(\sig(\bX)_{0t} =  \sig(\bX)_{t}\) and \(\sig^N(\bX)_{0t} = \sig^N(\bX)_{t}\).
\end{rem}

\begin{defn}[Shuffle product on \(T(\RR^d)\)]
Let \(e_1,\dots,e_d\) be the standard basis of \(\RR^d\). For a word \(\bw=w_1\dots w_n\), set \(e_{\bw}:=e_{w_1}\otimes\dots\otimes e_{w_n}\), with \(e_{\emptyset}:=1\).

For \(\mathbf{u}=\sum_{\bw\in\bW}u_{\bw}e_{\bw}\) and \(\mathbf{v}=\sum_{\bw'\in\bW}v_{\bw'}e_{\bw'}\) in \(T(\RR^d)\), define
\[\mathbf{u}\shuffle\mathbf{v}:=\sum_{\bw,\bw'\in\bW}u_{\bw}v_{\bw'}(e_{\bw}\shuffle e_{\bw'}),\]
where \(e_{\bw}\shuffle e_{\bw'}\) is obtained from the previously defined shuffle product of words \(\bw\shuffle\bw'\).
\end{defn}

\begin{defn}
    The space \(T(\RR^d)\) endowed with the shuffle product \(\shuffle\),  i.e., the quadruple \(\left(T(\RR^d), +, \cdot, \shuffle \right)\), is an associative algebra.
\end{defn}
The shuffle product encodes the multiplicative structure of the signature by expressing products of iterated integrals as linear combinations of iterated integrals. This yields the following fundamental identity for the signature of a path of bounded variation, known as the \emph{shuffle property}.

\begin{prop}[Shuffle Property] \cite{chevyrev2025primer} Let \((\bX_t)_{t \in [0,T]}\) be a continuous \(\RR^d\)-valued path of bounded variation and \(\bw = w_1 \dots w_n\) and \(\bw' = w'_1 \dots w'_m\) be two words then
\[\innerprod{\sig(\bX)_{st}, \bw} \innerprod{\sig(\bX)_{st}, \bw'} = \innerprod{\sig(\bX)_{st}, \bw \shuffle \bw'}, \quad \text{ for every } 0\le s\le t\le T.\]
\end{prop}
\begin{rem} \label{p_value_remark}
The shuffle property also holds for continuous paths of bounded \(p\)-variation with \(1\le p<2\), when the integrals are interpreted in Young's sense. It also holds for Brownian motion and fractional Brownian motion with Hurst parameter \(H>\frac{1}{4}\), when the iterated integrals are interpreted in the Stratonovich sense through their canonical geometric rough path lifts.
\end{rem}

We define the signature of a data stream as the signature of the piecewise linear path interpolating its data points. Computationally, we first compute the signature over each linear segment between consecutive data points. The signature over a longer interval is then obtained by combining these segment signatures using Chen's identity, which describes how signatures over adjacent intervals compose. We now state this identity formally.

\begin{prop}[Chen's Identity]
\cite{Friz_Victoir_2010}
Let \((\bX_t)_{t \in [0,T]}\) be a continuous, \(\mathbb{R}^d\)-valued path of bounded variation. Then, the concatenated signature over intervals \([s, u]\) and \([u, t]\) satisfies
\[\sig(\bX)_{st} = \sig(\bX)_{su} \otimes \sig(\bX)_{ut},\]
for each \(0 \leq s \leq u \leq t \leq T\). This identity can be equivalently expressed as follows:
\[
\innerprod{\sig(\bX)_{st}, \bw} = \sum_{\bw_1, \bw_2 \,:\, \bw_1 \bw_2 = \bw} \innerprod{\sig(\bX)_{su}, \bw_1} \innerprod{\sig(\bX)_{ut}, \bw_2} ,
\]
where \(\bw\) is an arbitrary word from \(\bW\).
\end{prop}

\subsection{Rooted trees and forests}

Let the set of all rooted trees (forests) with finitely many vertices be denoted by \(\trees\) (\(\mathcal{F}\)), and those with at most \(n\) vertices by \(\trees_n\) (\(\mathcal{F}_n\)). The vertices are decorated with letters from \(\cA=\{1,2,\dots,d\}\). The formal products of rooted trees form forests. The empty forest is denoted by \(\mathbf{1}\). The degree \(|h|\) is the number of vertices in a tree \(h\), and the degree of a forest is sum of the degrees of trees in it, with \(|\mathbf{1}|=0\).

Every tree can be constructed recursively as \([h_1\cdots h_n]_r\), by attaching \(h_1,\dots,h_n\) to a root decorated by \(r\). For example,
\[[\mathbf{1}]_i=\Forest{Decorated[i]},\qquad [\Forest{Decorated[i]}]_j=\Forest{Decorated[j[i]]},\qquad [\Forest{Decorated[i]}\Forest{Decorated[j]}]_k=\Forest{Decorated[k[i][j]]}.\]
We assume that the order of the branches does not matter, that is, \([h_1\cdots h_n]_r=[h_{\sigma(1)}\cdots h_{\sigma(n)}]_r\) for any permutation \(\sigma\).

\subsection{The Connes--Kreimer Hopf algebra}
The Connes--Kreimer Hopf algebra \((\hopf,\cdot,\Delta,\cS)\) is the commutative polynomial algebra generated by the trees in \(\trees\). Its basis consists of forests, and its product is forest multiplication, with unit \(\mathbf{1}\). This algebra is used in \cite{connes1999hopf} in the context of renormalization theory.

The coproduct is defined recursively by \(\Delta\mathbf{1}=\mathbf{1}\hat{\otimes}\mathbf{1}\) and
\begin{equation}\label{coproduct h}
\Delta[h_1\cdots h_n]_r=[h_1\cdots h_n]_r\hat{\otimes}\mathbf{1}+\big(\Id\hat{\otimes}[\,\cdot\,]_r\big)\Delta(h_1\cdots h_n).
\end{equation}
It preserves forest multiplication, \(\Delta(h_1\cdots h_n)=\Delta h_1\cdots\Delta h_n\), and is coassociative, \((\Delta\hat{\otimes}\Id)\Delta=(\Id\hat{\otimes}\Delta)\Delta\).

The counit, denoted by \(\mathbf{1}^*\), and the antipode \(\cS\) satisfy the following identity
\[P\big((\cS\hat{\otimes}\Id)\Delta h\big)=P\big((\Id\hat{\otimes}\cS)\Delta h\big)=\innerprod{\mathbf{1}^*,h}\mathbf{1},\]
where \(P(h_1\hat{\otimes}h_2)=h_1h_2\) and \(\Id\) is the identity map. The grading \(\hopf=\bigoplus_{n\ge0}\hopf_{(n)}\) is given by the number of vertices. Since \(\hopf_{(0)}=\RR\mathbf{1}\), this graded bialgebra is a Hopf algebra with a unique antipode.

The coproduct induces the convolution product \(\star\) on \(\hopf^*\), defined by
\[\innerprod{f\star g,h}=\innerprod{f\hat{\otimes}g,\Delta h}=\sum_{(h)}\innerprod{f,h^{(1)}}\innerprod{g,h^{(2)}},\]
where \(\Delta h=\sum_{(h)}h^{(1)}\hat{\otimes}h^{(2)}\). For further details, we refer to \cite{abe2004hopf,brouder2004trees,manchon2008hopf}.

\subsection{The branched signature}
For non-geometric rough paths, the usual shuffle property need not hold, so we retain products of iterated integrals explicitly. We use rooted trees instead of words to encode these iterated integrals and their products. Let \(\bX:[0,T]\to\RR^d\) be an \(\alpha\)-H\"older continuous path for \(\alpha > \frac{1}{4}\). The iterated integrals and their products below are formal expressions and can be made precise once a suitable notion of integration has been specified. Following the recursive construction in \cite{gubinelli2010ramification}, we now define the branched signature by collecting these iterated integrals and their products into a functional on the Connes--Kreimer Hopf algebra.

\begin{defn}[Branched signature]
Let \(\cA=\{1,\dots,d\}\) be the alphabet of decorations for a given \(d\)-dimensional path \(\bX:[0,T]\to\RR^d\). Denote by \(\trees\) the set of rooted trees with vertices decorated in \(\cA\), and by \(\mathcal{F}\) the corresponding set of forests. Let \(\hopf\) be the decorated Connes--Kreimer Hopf algebra generated by \(\trees\), with product given by disjoint union of forests and unit \(\mathbf{1}\).

The \emph{branched signature} of \(\bX\), denoted by \(\bsig(\bX)_{st}\), is a functional on \(\hopf\) and is defined as follows
\begin{equation}
\bsig(\bX)_{st}=\sum_{\tau\, \in\, \mathcal{F},\,|\tau|\, \ge \, 0}\innerprod{\bsig(\bX)_{st},\tau}\se_{\tau},
\end{equation}
where \(\se_{\tau}\) is the coordinate functional dual to the forest \(\tau\). The iterated integrals below are understood formally until a suitable notion of integration has been specified.

For \(0\le s\le t\le T\), the coordinates are recursively defined by
\[\innerprod{\bsig(\bX)_{st},\mathbf{1}}=1,\qquad \innerprod{\bsig(\bX)_{st},\tau}=\int_s^t\innerprod{\bsig(\bX)_{su},\tau'}\,d\bX_u^{\mathbf{r}},\]
where \(\tau\) is a nonempty rooted tree, \(\mathbf{r}\) is the decoration of its root, and \(\tau'\) is the forest obtained by removing the root from \(\tau\). For a forest \(\tau_1\cdots\tau_m\), we set
\[\innerprod{\bsig(\bX)_{st},\tau_1\cdots\tau_m}=\prod_{j=1}^m\innerprod{\bsig(\bX)_{st},\tau_j}.\]
These coordinates are extended linearly to \(\hopf\).
\end{defn}

\begin{rem}
Under the assumptions of the extension theorem, a truncated branched rough path admits a unique extension to all degrees satisfying the required regularity bounds, which defines its branched signature \cite{gubinelli2010ramification}. This is analogous to the unique extension of a geometric rough path to its signature \cite{Friz_Victoir_2010}.
\end{rem}

For a suitable choice of integration, the corresponding Chen's identity is expressed through the coproduct as follows.

\begin{prop}[Chen's identity for the branched signature]
For \(0\le s\le u\le t\le T\) and every forest \(h\in\hopf\),
\[\innerprod{\bsig(\bX)_{st},h}=\sum_{(h)}\innerprod{\bsig(\bX)_{su},h^{(1)}}\innerprod{\bsig(\bX)_{ut},h^{(2)}}.\]
Equivalently, \(\bsig(\bX)_{st}=\bsig(\bX)_{su}\star\bsig(\bX)_{ut}\).
\end{prop}

\subsection{The classical Volterra signature}
Following \cite{harang2021volterra,hager2026volterra}, we introduce iterated integrals weighted by a Volterra kernel, and hence the classical Volterra signature. We begin with a continuously differentiable path \(\bX:[0,T]\to\RR^d\) and a measurable scalar kernel \(k(\tau,r)\), defined for \(0\le r<\tau\le T\), such that
\[\sup_{\tau \, \in \, [0,T]}\int_0^\tau |k(\tau,r)|\,dr<\infty.\]
The recent paper on the Volterra signature \cite{hager2026volterra} talks about matrix-valued kernels. To keep things simple we restrict ourselves to the scalar kernels only. Throughout this subsection, the kernel \(k\) is fixed. Since \(\bX\) is continuously differentiable, we write \(d\bX_r=\dot{\bX}_r\,dr\). For \(0\le s\le t\le\tau\le T\), the first-level Volterra integral is defined by
\[z^\tau_{st}:=\int_s^t k(\tau,r)\,d\bX_r.\]
Here \(\tau\) is an external time parameter, which is kept fixed while integrating over \([s,t]\).

\begin{defn}[The classical Volterra signature]
The classical Volterra signature \(\vsig(\bX)^\tau_{st}\) of a continuously differentiable path \(\bX:[0,T] \to \RR^d\) is a \(T((\RR^d))\)-valued process and is defined through its coordinates \(\vsig(\bX)^{\bw,\tau}_{st}:=\innerprod{\vsig(\bX)^\tau_{st},\bw}\) for a word \(\bw\). We set \(\vsig(\bX)^{\emptyset,\tau}_{st}=1\), and for a word \(\bw\) and a letter \(i\in\cA\), define
\[\vsig(\bX)^{\bw i,\tau}_{st}:=\int_s^t k(\tau,r)\vsig(\bX)^{\bw,r}_{sr}\,d\bX_r^i.\]
In particular, \(\vsig(\bX)^{i,\tau}_{st}\) is the \(i\)-th coordinate of \(z^\tau_{st}\).
\end{defn}
To make sense of the Chen's identity in the case of Volterra signature, we use the convolution product, which is defined below.
\begin{defn}[The convolution product]
Let \(z\) be the first-level Volterra integral defined above. For a bounded measurable function \(y:[s,T]\to\RR\), set \([y\ast z]^{\emptyset,\tau}_{st}:=y_\tau\) and define recursively
\[[y\ast z]^{\bw i,\tau}_{st}:=\int_s^t k(\tau,r)[y\ast z]^{\bw,r}_{sr}\,d\bX_r^i.\]
Thus, the external parameter in the earlier factor is evaluated at an integration time. Taking \(y\equiv1\) gives \([1\ast z]^{\bw,\tau}_{st}=\vsig(\bX)^{\bw,\tau}_{st}\).
\end{defn}
This convolution product replaces ordinary tensor multiplication in the classical Chen's identity.
\begin{prop}[The Volterra Chen's identity]
For \(0\le s\le u\le t\le\tau\le T\) and a word \(\bw=w_1\cdots w_n\),
\[\vsig(\bX)^{\bw,\tau}_{st}=\sum_{k=0}^n\left[\vsig(\bX)^{w_1\cdots w_k,\cdot}_{su}\ast z\right]^{w_{k+1}\cdots w_n,\tau}_{ut}.\]
Here the first factor is the function \( \lambda \mapsto \vsig(\bX)^{w_1\cdots w_k,\lambda}_{su}\), and empty prefixes and suffixes are interpreted as \(\emptyset\).
\end{prop}

\begin{rem}
If \(k\equiv1\), the external time parameter vanishes and we get \(\vsig(\bX)^\tau_{st}=\sig(\bX)_{st}\).
\end{rem}

For irregular driving paths, the first-level Volterra integral may not determine the higher iterated integrals. A suitable Volterra rough path lift provides the additional information needed to define these integrals and the corresponding convolution.

\section{The generalized Volterra signature}
In this section, we are going to introduce the generalized Volterra signature and the associated properties. Everything in this section should be considered formal. Any technical results about the existence of the Volterra signature are postponed to the succeeding sections. To begin with, we introduce some terminology before defining the generalized Volterra signature formally. 
\begin{defn}[Kernel tree]
The kernel tree \(\ktree\) is the tree that determines the arguments in the product of the kernels in the Volterra signature. Always a decorated tree with no or fixed decoration on the fixed root for the parameter \(\tau\) in the kernel. Some examples of lower degree kernel trees are
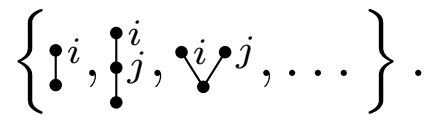
\[
\left\{\Forest{Decorated[[i]]}, \Forest{Decorated[[j[i]]]}, \Forest{Decorated[[i][j]]}, \dots  \right\}.
\]
For instance, corresponding to the kernel tree \(\Forest{Decorated[[j[i]]]}\) the kernel product is
\[
k(\tau, r_j)k(r_j, r_i).\]
\end{defn}

\begin{defn}[Differential tree]
The differential tree \(\dtree\) is the tree that determines the order of the differentials and the domain of the integration. Always a decorated rooted tree. Some examples of differential trees are
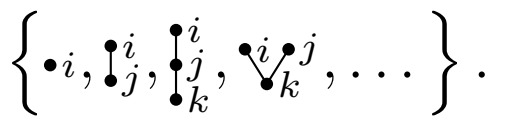
\[
\left\{\Forest{Decorated[i]}, \Forest{Decorated[j[i]]}, \Forest{Decorated[k[j[i]]]}, \Forest{Decorated[k[i][j]]}, \dots \right\}.
\]
\end{defn}

\begin{defn}[Volterra tree]
The Volterra tree \(\vtree\) is the pair \((\ktree, \dtree)\) that contains the kernel tree \(\ktree\) and the differential tree \(\dtree\).
\end{defn}
Every pair of the kernel and the differential trees is not useful. We stick to the one that gives us a well--defined notion of the integration in terms of the limit and the relationship between the kernel and the components of the path we are integrating against. Therefore, we only work with compatible Volterra trees defined below.
\begin{defn}[The compatible Volterra tree] A Volterra tree \(\vtree :=(\ktree, \dtree)\) is compatible if the following conditions hold
\begin{enumerate}
    \item \(N_{\ktree^*} = N_{\dtree}\), where \(N_{\ktree^*}\) and \(N_{\dtree}\) are the sets of vertices in the trees \(\ktree^*\) and \(\dtree\) respectively and \(\ktree^*\) is the tree after removing the dummy root from the kernel tree \(\ktree\).
    \item If \(v = \kparent{u} \ne z\) then \(u \prec_{\dtree} v\), where  \(\kparent{u}\) is the parent of vertex \(u\) in the kernel tree and is not the root vertex, and \(\prec_{\dtree}\) is the partial order determined by the differential tree \(\dtree\).
\end{enumerate}
\end{defn}

\begin{example}
The Volterra tree
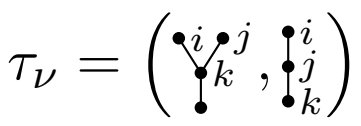
\[\vtree = \left(\Forest{Decorated[[k[i][j]]]}, \Forest{Decorated[k[j[i]]]} \right) \]
is compatible. Indeed, in the kernel tree, \(\kparent{i}=\kparent{j}=k, \) while the partial order induced by the differential tree satisfies
\(i\prec_{\dtree}k\) and\(j\prec_{\dtree}k.\)
\end{example}

\begin{example}[An incompatible Volterra tree]
The Volterra tree
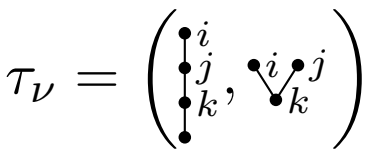
\[\vtree = \left(\Forest{Decorated[[k[j[i]]]]}, \Forest{Decorated[k[i][j]]} \right) \]
is incompatible because in the kernel tree \(j = \kparent{i}\) but in the differential tree there is no relationship between the vertices \(i\) and \(j\) with respect to partial order \(\prec_{\dtree}\) determined by the differenital tree.
\end{example}

\begin{defn}[The domain of integration]
    The domain of integration \(\bA_{st}^{\dtree}\) determined by the differential tree \(\dtree\) is defined as
    \[
\bA_{st}^{\dtree} := \left\{ (r_v)_{v \in N_{\dtree}}: r_u < r_v \text{ whenever } u \prec_{\dtree} v \text{ and } s < r_v < t \text{ for all } v \in N_{\dtree} \right\}.
\]
\end{defn}

\begin{defn}[Kernel function] \label{kernel_function}
Let the \(2\)--simplex be
\[ \Delta_T^2 := \left\{(s,t)\in[0,T]^2:0\leq s \leq t\leq T\right\}.\]
A function \[ k:\Delta_T^2\to\RR \] is called a \emph{kernel function}. The kernel \(k\) is allowed to be singular along the diagonal
\[ \{(s,t)\in[0,T]^2:s=t\}. \] 
However, we will stick to continuous or at least integrable kernels. Unless otherwise stated, we stick to only scalar--valued kernels. The matrix valued kernels are studied in the classical Volterra signature setting in \cite{hager2026volterra}.
\end{defn}
Next, we are going to establish the algebraic structure on the compatible Volterra trees and forests that is required to define generalized Volterra signature and the corresponding identities. In particular, we are going to define Hopf algebra of compatible Volterra forests, with the corresponding operations, i.e., coproduct, counit, and antipode.

\subsection{Hopf algebra of compatible Volterra forests}
In this subsection we are going to introduce a natural algebraic structure that encodes generalized Volterra signature \(\gvsig\). With this in hand, \(\gvsig\) will be a character on an algebra generated by compatible Volterra forests, an extension of compatible Volterra trees, that we are going to introduce in this section. 

Let us introduce more terminology after we introduced the the kernel tree \(\ktree\), the differential tree \(\dtree\), and the Volterra tree \(\vtree\). 
\begin{defn}[Compatible Volterra forest]
    Let the kernel forest and the the differential forest be denoted by \(\kforest\) and \(\dforest\) respectively, where the forest is simply formal product of the trees. We define the compatible Volterra forest \(\vforest\) to be the pair \((\kforest, \dforest)\). Alternatively, if \(\{\vtree^i = (\ktree^i, \dtree^i)\}_{i=1}^{m}\) is a collection of \(m\) compatible Volterra trees. Then the forest made of these trees is 
    \[
    \vforest = \vtree^1 \cdots \vtree^m = (\ktree^i, \dtree^i) \cdots (\ktree^i, \dtree^i).
    \]
    This is just the formal notation. In kernel and differential tree notation the forest is written as
    \[
    \vforest = \left(\ktree^1 \dots \ktree^m, \dtree^1 \dots \dtree^m\right),
    \]
    where the kernel and the differential tree product is usual forest product. The dummy root in the kernel trees can be identified to a single vertex too.
\end{defn}

Let \(\KK\) be the field of characteristic zero, say \(\KK= \RR\). Let the set of compatible Volterra trees and forests be denoted by \(\tset\) and \(\fset\) respectively. Let the empty forest be denoted by \(\mathbf{1}\). We define the degree of a compatible Volterra tree as
\[
\deg(\vtree) = |N_{\dtree}| = |N_{\ktree^*}|.
\]
Whereas, the degree of a compatible Volterra forest \(\vforest = \vtree^1 \dots \vtree^k\) is defined as 
\[
\deg(\vforest) = \deg(\vtree^1) + \dots + \deg(\vtree^k).
\]
The vector space spanned by the set of compatible Volterra forest is defined as
\[
\hopf_{\nu} := \myspan_{\RR}(\fset),
\]
where the underlying filed is \(\RR\). The space \(\hopf_{\nu}\) is graded in the sense that 
\[
\hopf_{\nu} = \bigoplus_{n \ge 0} \hopf_{\nu}^n,
\]
where 
\[
\hopf_{\nu}^n := \myspan_{\RR}\{\vforest \in \fset: \deg(\vforest) = n\} \quad \text{ and } \quad \hopf_{\nu}^0 := \RR.\mathbf{1}.
\]
The space \(\hopf_{\nu}\) is equipped the forest product defined as
\[
\mathbf{1}\cdot\vforest = \vforest\cdot\mathbf{1} = \vforest \quad \text{ and } \quad  \vforest^1 \cdot \vforest^2 := \left( \kforest^1 \sqcup\kforest^2,  \dforest^1\dforest^2\right),
\]
where \(\kforest^1\sqcup\kforest^2\) denotes the kernel forest obtained by merging the two dummy roots into a single root and \(\dforest^1\dforest^2\) denotes the disjoint union of the differential forests \(\dforest^1\) and \(\dforest^2\).

The forest product is commutative as the trees we are dealing with are non-planar i.e.,
\[
\vforest^1 \cdot \vforest^2 = \vforest^2 \cdot \vforest^1.
\]
It is also associative i.e., \(\vforest^1 \cdot \left(\vforest^2 \cdot \vforest^3\right) = \left( \vforest^1 \cdot \vforest^2 \right) \cdot \vforest^3\). Also, 
\[
\deg \left(\vforest^1 \cdot \vforest^2 \right) =\deg \left( \vforest^1\right) +  \deg \left(\vforest^2\right).
\]

\begin{defn}
    The space \(\left( \hopf_{\nu}, \cdot, \mathbf{1}\right)\), equipped with the forest product \(\cdot\) and the unit \(\mathbf{1}\),  is a graded commutative unital algebra.
\end{defn}

Next step in our discussion is to define the Hopf algebra of the compatible Volterra forests. To do this, we need to define some operations on the compatible Volterra forests. To begin with, let us define the set of admissible cuts and the coproduct. 

\begin{defn}[Set of admissible cuts]
    A cut \(c\) of a differential tree \(\dtree\) is a collection of all the edges such that a path from the root to the leaf contains at most one edge of the cut \(c\). Collection of all the cuts of a differential tree \(\dtree\) is denoted by \(\operatorname{Adm}(\dtree)\).
\end{defn}
\begin{example}
    For the differential tree \(\dtree = \Forest{Decorated[k[i][j]]}\), the  set of admissible cuts comprise the following.
    \[
    c_1 = \{(i,k)\}, \quad c_2 = \{(j,k)\}, \quad \text{ and } \quad c_3 = \{(i,k),(j,k)\}.
    \]
\end{example}

\(\operatorname{Adm}(\dtree)\) is also called the set of nontrivial admissible cuts of the differential tree \(\dtree\). A cut \(c\) divides the tree into two parts called pruned forest \(P^c(\dtree)\) and the trunk \(R^c(\dtree)\). The pruned forest \(P^c(\dtree)\) contains the set of vertices lying strictly before the cut in the sense of partial order determined by the tree. The trunk \(R^c(\dtree)\) is the remaining part that contains the root of the tree.

\begin{example}
    For the differential tree \(\dtree = \Forest{Decorated[l[k[i][j]]]}\) and for one of its cuts \(c = \{(i,k),(j,k)\}\), the \(P^c(\dtree)\) and \(R^c(\dtree)\) are given as
    \[
    P^c(\dtree) = \Forest{Decorated[i]} \, \Forest{Decorated[j]} \quad \text{ and } \quad R^c(\dtree) = \Forest{Decorated[l[k]]}.
    \]
\end{example}

The same cut \(c\), the differential pruned forest \(P^c(\dtree)\), and the differential trunk \(R^c(\dtree)\) of the tree \(\dtree\) also determine the kernel pruned forest and the kernel trunk of the corresponding kernel tree \(\ktree\). The kernel pruned forest \(P^c(\ktree)\) and the kernel trunk \(R^c(\ktree)\) is simply the restriction of the kernel tree \(\ktree\) to the differential pruned forest \(P^c(\dtree)\) and the differential trunk \(R^c(\dtree)\) respectively with an additional dummy root in each tree of the kernel pruned forest \(P^c(\ktree)\). We retain the actual dummy root of the kernel tree \(\ktree\) in the kernel trunk \(R^c(\ktree)\).

\begin{example}
    For the Volterra tree \(\vtree = \left(\Forest{Decorated[[k[i][j]]]}, \Forest{Decorated[k[j[i]]]}\right)\) and the cut \(c = \{(j,k)\}\), we have
    \[
    P^c(\ktree) = \Forest{Decorated[[i]]} \, \Forest{Decorated[[j]]} \quad \text{ and } \quad R^c(\ktree) = \Forest{Decorated[[k]]}.
    \]
\end{example}


\begin{defn}[Pruned forest and trunk of the compatible Volterra trees]
For a compatible Volterra tree \(\vtree = (\ktree, \dtree)\) and an admissible cut \(c\), we define
\[
P^c(\vtree):=\left(P^c(\ktree),P^c(\dtree)\right)
\]
and
\[
R^c(\vtree):=\left(R^c(\ktree),R^c(\dtree)\right).
\]
\end{defn}
\begin{defn}[Coproduct of the compatible Volterra tree \(\vtree\)]
For a compatible Volterra tree \(\vtree\), we define the coproduct \(\Delta: \hopf_{\nu} \to \hopf_{\nu} \otimes \hopf_{\nu}\) as
\[
\Delta\vtree:=\vtree\otimes\mathbf{1} + \mathbf{1}\otimes\vtree + \sum_{c \, \in \, \operatorname{Adm}(\vtree)} R^c(\vtree)\otimes P^c(\vtree),
\]
where \(\operatorname{Adm}(\vtree)\) denotes the set of nontrivial admissible cuts of the differential tree of \(\vtree\).
\end{defn}

The definition can be extended to forests multiplicatively as 
\[
\Delta(\vtree^1 \,  \vtree^2) =  \Delta(\vtree^1) \cdot \Delta(\vtree^2),
\]
where the product on the right hand side is defined as
\[
\left( \vtree^{i_1} \otimes \vtree^{i_2}\right) \left( \vtree^{j_1} \otimes \vtree^{j_2}\right) = \vtree^{i_1} \cdot \vtree^{j_1} \otimes \vtree^{i_2} \cdot \vtree^{j_2}, \quad \forall \, \vtree^{i_1},\vtree^{i_2},\vtree^{j_1},\vtree^{j_2} \in \tset.
\]
The coproduct is also coassociative in the following sense
\[
\left(\Delta \otimes \operatorname{Id}\right) \Delta = \left(\operatorname{Id} \otimes \Delta\right) \Delta.
\]

\begin{example}
Here are a few examples of the coproduct of the compatible Volterra trees of lower degree.
    \begin{enumerate}
        \item For \(\vtree = \left(\Forest{Decorated[[j[i]]]}, \Forest{Decorated[j[i]]}\right)\), we have
        \[\Delta \vtree = 1 \otimes  \vtree +  \vtree \otimes 1 + \left(\Forest{Decorated[[j]]}, \Forest{Decorated[j]}\right) \otimes \left(\Forest{Decorated[[i]]}, \Forest{Decorated[i]}\right).\]
        \item For \(\vtree = \left(\Forest{Decorated[[j][i]]}, \Forest{Decorated[j[i]]}\right)\), we have
        \[\Delta \vtree = 1 \otimes  \vtree +  \vtree \otimes 1 + \left(\Forest{Decorated[[j]]}, \Forest{Decorated[j]}\right) \otimes \left(\Forest{Decorated[[i]]}, \Forest{Decorated[i]}\right).\]
        \item For \(\vtree = \left(\Forest{Decorated[[k[j[i]]]]}, \Forest{Decorated[k[j[i]]]}\right)\), we have
        \[\Delta \vtree = 1 \otimes  \vtree +  \vtree \otimes 1 + \left(\Forest{Decorated[[k[j]]]}, \Forest{Decorated[k[j]]}\right) \otimes \left(\Forest{Decorated[[i]]}, \Forest{Decorated[i]}\right) + \left(\Forest{Decorated[[k]]}, \Forest{Decorated[k]}\right) \otimes \left(\Forest{Decorated[[j[i]]]}, \Forest{Decorated[j[i]]}\right).\]
        \item For \(\vtree = \left(\Forest{Decorated[[k[i][j]]]}, \Forest{Decorated[k[i][j]]}\right)\), we have
        \[\Delta \vtree = 1 \otimes  \vtree +  \vtree \otimes 1 + \left(\Forest{Decorated[[k[j]]]}, \Forest{Decorated[k[j]]}\right) \otimes \left(\Forest{Decorated[[i]]}, \Forest{Decorated[i]}\right) + \left(\Forest{Decorated[[k[i]]]}, \Forest{Decorated[k[i]]}\right) \otimes \left(\Forest{Decorated[[j]]}, \Forest{Decorated[j]}\right) + \left(\Forest{Decorated[[k]]}, \Forest{Decorated[k]}\right) \otimes \left(\Forest{Decorated[[i][j]]}, \Forest{Decorated[i]}\Forest{Decorated[j]}\right).\]
    \end{enumerate}
\end{example}

\begin{defn}[Counit]
    The counit \(\varepsilon\) is a map \(\varepsilon: \hopf_{\nu} \to \RR\) defined as
    \[
    \varepsilon(\vforest) = 
    \begin{cases}
        1, & \text{ if } \vforest = \mathbf{1},\\
        0, & \text{ otherwise.}
    \end{cases}
    \]
    The counit \(\varepsilon\) satisfies the following identity
\[
\left(\varepsilon \otimes \operatorname{Id}\right) \Delta = \operatorname{Id} = \left(\operatorname{Id} \otimes \, \varepsilon\right) \Delta.
\]
\end{defn}

\begin{defn}[Antipode]
    The antipode \(\cS\) is a map \(\cS: \hopf_{\nu} \to\hopf_{\nu} \) that encodes the inverses and is defined recursively as
    \[
    \cS(\mathbf{1})=\mathbf{1}, \quad \cS(\vtree) = - \, \vtree \,  - \sum_{c \, \in \, \operatorname{Adm}(\vtree)} R^c(\vtree) \, \cS \left(P^c(\vtree)\right).
    \]
The antipode can also be extended to forests multiplicatively as follows.
\[
\cS (\vtree^1 \, \vtree^2) = \cS (\vtree^1)  \, \cS (\vtree^2), \quad \forall \,  \vtree^1, \vtree^2 \in \tset.
\]
Moreover, the following diagram commutes.
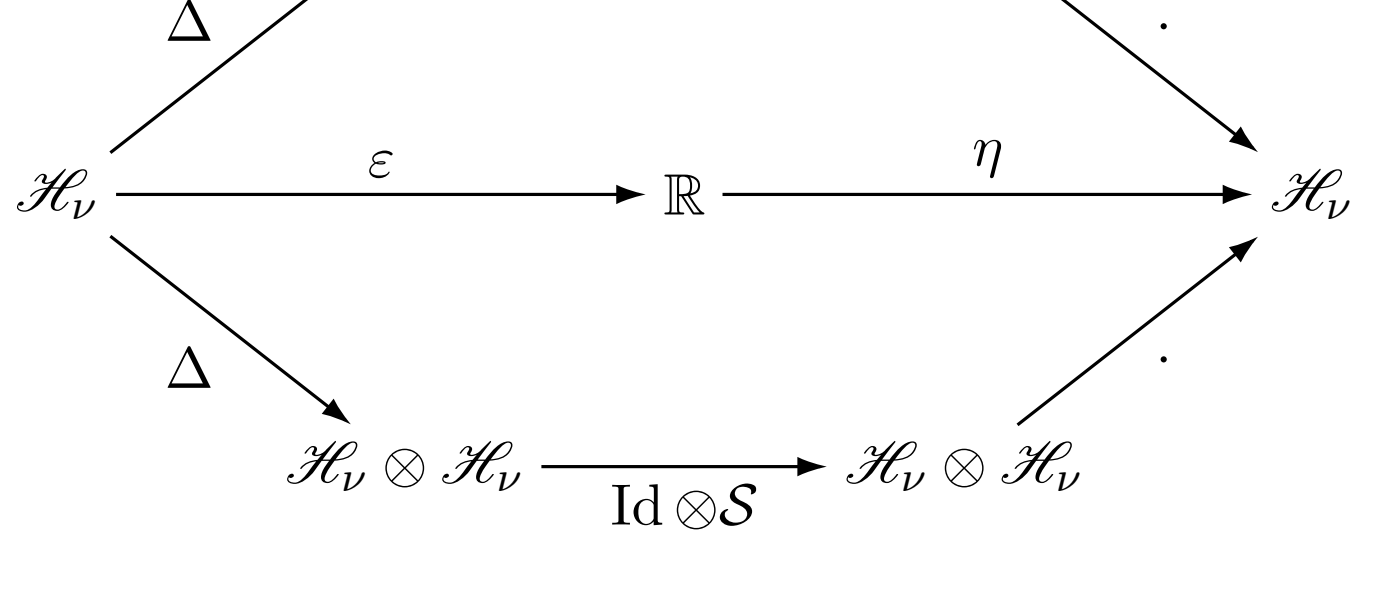
\begin{center}
\begin{tikzpicture}[
    >={Latex[length=2mm,width=1.3mm]},
    every node/.style={inner sep=3pt},
    every edge/.style={draw,->,semithick}
]
    \node (left)        at (0,0)       {$\hopf_{\nu}$};
    \node (upperleft)   at (2.3,1.8)   {$\hopf_{\nu}\otimes\hopf_{\nu}$};
    \node (upperright)  at (6.0,1.8)   {$\hopf_{\nu}\otimes\hopf_{\nu}$};
    \node (field)       at (4.15,0)    {$\RR$};
    \node (right)       at (8.3,0)     {$\hopf_{\nu}$};
    \node (lowerleft)   at (2.3,-1.8)  {$\hopf_{\nu}\otimes\hopf_{\nu}$};
    \node (lowerright)  at (6.0,-1.8)  {$\hopf_{\nu}\otimes\hopf_{\nu}$};

    \path
    (left)       edge node[above left]  {$\Delta$} (upperleft)
    (upperleft)  edge node[above]       {$\cS\otimes\operatorname{Id}$} (upperright)
    (upperright) edge node[above right] {$\cdot$} (right)
    (left)       edge node[above]       {$\varepsilon$} (field)
    (field)      edge node[above]       {$\eta$} (right)
    (left)       edge node[below left]  {$\Delta$} (lowerleft)
    (lowerleft)  edge node[below]       {$\operatorname{Id}\otimes\cS$} (lowerright)
    (lowerright) edge node[below right] {$\cdot$} (right);
\end{tikzpicture}
\end{center}
Here $\eta:\RR\to\hopf_{\nu}$ is the unit map,
defined by $\eta(a)=a\mathbf{1}$, and the arrows
labelled $\cdot$ denote the multiplication map
\[
\cdot:\hopf_{\nu}\otimes\hopf_{\nu}\to\hopf_{\nu}.
\]
The diagram commutes in the sense that
\[
\cdot\circ(\cS\otimes\operatorname{Id})\circ\Delta
=
\eta\circ\varepsilon
=
\cdot\circ(\operatorname{Id}\otimes\cS)\circ\Delta.
\]
\end{defn}




\begin{defn}[Hopf algebra of compatible Volterra forests]
    The space \(\left(\hopf_{\nu}, \cdot, \mathbf{1}, \Delta, \varepsilon, \cA \right)\) is a connected graded commutative Hopf algebra. 
\end{defn}

\subsection{The generalized Volterra signature}
Let us define the notion of generalized Volterra signature after we have a enough terminology built on compatible Volterra trees and the corresponding algebraic structure. For now, we will define Volterra signature for the paths of bounded variation where the integrals are defined in the Riemann-Steiltjes sense. The extension of the definition to other paths is postponed to our later works.

\begin{defn}\label{gen_volterra_component}
Let \(\bX : [0,T] \to \RR^d\) be the path of bounded variation. Let \(\vtree :=(\ktree, \dtree)\) be a compatible Volterra tree. For \(0 \le s \le t \le \tau \le T\), we define the component of the generalized Volterra signature \(\gvsig(\bX)^{\tau}_{st}\) of \(\bX\) corresponding to the tree \(\vtree\) as follows.
\begin{equation}
\innerprod{\gvsig(\bX)^{\tau}_{st}, \vtree} := \int_{\bA_{st}^{\dtree}} \left(\prod_{v \in N_{\ktree^*}} k (r_{\kparent{v}}, r_v)\right) \prod_{v \in N_{\dtree}} d\bX_{r_v}^{v},
\end{equation}
where \(k(s,t)\) is the kernel function defined in \ref{kernel_function}, \(\kparent{v}\) is the parent of vertex \(v\) in the kernel tree \(\ktree\), and \(\ktree^*\) is the kernel tree without the dummy root.
\end{defn}
Furthermore, for the empty forest \(\mathbf{1}\), we set \( \innerprod{\gvsig(\bX)_{st}^{\tau},\mathbf{1}}:=1 \). Also,
the definition extends naturally from Volterra trees to Volterra forests. Let \(\vtree = \vtree^1 \dots \vtree^k\), where \(\vtree^l = (\ktree^l, \dtree^l)\), \(l=1,2,\dots, k\), be a Volterra forest whose components \( \vtree^1,\dots,\vtree^k\) are compatible Volterra trees. We define its generalized Volterra signature component multiplicatively as follows.
\[
\innerprod{\gvsig(\bX)^{\tau}_{st},\vtree} := \prod_{l=1}^{k}
\innerprod{\gvsig(\bX)^{\tau}_{st},\vtree^l}.
\]
\begin{example}
For \(\ktree = \Forest{Decorated[[k[i][j]]]}\) and \(\dtree = \Forest{Decorated[k[j[i]]]}\) i.e., \(\vtree = \left(\Forest{Decorated[[k[i][j]]]},\Forest{Decorated[k[j[i]]]} \right)\), we have
\[
\bA_{st}^{\dtree} = \{(r_i, r_j, r_k) \in [s,t]^3: s \le r_i  \le r_j \le r_k \le t\}.
\]
The component of the \(\gvsig(\bX)^{\tau}_{st}\) corresponding to \(\vtree\) is given below
\[
\innerprod{\gvsig(\bX)^{\tau}_{st}, \vtree} = \int_{s \le r_i  \le r_j \le r_k \le t} k(\tau, r_k) k(r_k,r_j) k(r_k,r_i) d\bX^{\Forest{Decorated[i]}}_{r_i} d\bX^{\Forest{Decorated[j]}}_{r_j} d\bX^{\Forest{Decorated[k]}}_{r_k}.
\]
Equivalently,
\[
\innerprod{\gvsig(\bX)^{\tau}_{st}, \vtree} = \int_{s}^{t} k(\tau, r_k) \left(\int_{s}^{r_k} k(r_k,r_j) d\bX^{\Forest{Decorated[j]}}_{r_j}\right) \left(\int_{s}^{r_j} k(r_k,r_i) d\bX^{\Forest{Decorated[i]}}_{r_i}\right) d\bX^{\Forest{Decorated[k]}}_{r_k}.
\]
\end{example}
\begin{example}[An incompatible Volterra tree]
Consider
\[
\vtree=
\left(
\Forest{Decorated[[i[j[k]]]]},
\Forest{Decorated[k[j[i]]]}
\right).
\]
The differential tree gives
\[
\bA_{st}^{\dtree}
=
\left\{
(r_i,r_j,r_k)\in[s,t]^3:
s\leq r_i\leq r_j\leq r_k\leq t
\right\}.
\]
The formal expression associated with this pair is
\[
\begin{aligned}
\innerprod{\gvsig(\bX)^{\tau}_{st},\vtree}
& = \int_s^t k(\tau,r_i)\,
\left[  \int_s^{r_k} k(r_i,r_j)\,
\left(\int_s^{r_j} k(r_j,r_k)\,
d\bX^{\Forest{Decorated[i]}}_{r_i}
\right) d\bX^{\Forest{Decorated[j]}}_{r_j}
\right]
d\bX^{\Forest{Decorated[k]}}_{r_k}.
\end{aligned}
\]
This is not a valid generalized Volterra signature component because the integration order determined by the differential tree \( \dtree = \Forest{Decorated[k[j[i]]]}\) is inconsistent with the parent--child order determined by the kernel tree \( \ktree  =\Forest{Decorated[[i[j[k]]]]}\).
\end{example}

\begin{defn}[The generalized Volterra signature]\label{gen_volterra_def}
Let \(\cT_{\nu}\) denote the set of all compatible Volterra trees, and let \(\cF_{\nu}\) denote the set of all compatible Volterra forests, including the empty forest \(\mathbf{1}\). The generalized Volterra signature of \(\bX\) over \([s,t]\) is the formal series
\[ \gvsig(\bX)_{st}^{\tau} = \mathbf{1} + \sum_{\substack{\vtree\in\mathcal{F}_{\nu}\\ |\vtree|\geq 1 }}
\innerprod{ \gvsig(\bX)_{st}^{\tau}, \vtree } \vtree, \]
where every forest has the form \( \vtree = \vtree^1\cdots\vtree^m \) with \(|\vtree|:= \sum_{j=1}^{m}
\left|N_{\dtree^j}\right|.\)
Also, we define the generalized Volterra signature of \(\bX\), truncated at step--\(N\), by the finite series
\[ \gvsig^N(\bX)_{st}^{\tau} = \mathbf{1} + \sum_{\substack{\vtree \, \in \, \mathcal{F}_{\nu}\\ 1 \, \leq \, |\vtree| \, \leq \, N }}
\innerprod{ \gvsig(\bX)_{st}^{\tau}, \vtree } \vtree.\]
\end{defn}
\noindent
We also define the generalized Volterra signature model \(\bV_{\ell}^{N}(\bX)_{st}^{\tau}\) as follows.
\begin{defn}[Generalized Volterra signature model]
Let \(\hopf_{\nu}\) be the Hopf algebra of compatible Volterra forests, and let \(\fset\) denote the set of these forests, including the empty forest \(\mathbf{1}\). Denote the subspace \(\hopf^{N}_{\nu}:=\operatorname{span}\{\vtree\in \fset:|\vtree|\leq N\}\). Then, for coefficients \(\ell=\ell_{\mathbf{1}}\mathbf{1}+\sum_{1 \, \leq \, |\vtree| \, \leq \, N}\ell_{\vtree}\vtree\in \hopf^{N}_{\nu} \), the generalized Volterra signature model \(\bV_{\ell}^{N}(\bX)_{st}^{\tau}\) (truncated to level \(N\)) for a \(d\)-dimensional path \(\bX:[0,T]\to\mathbb{R}^{d}\) is the linear functional of the generalized Volterra signature, i.e.,
\[
\bV_{\ell}^{N}(\bX)_{st}^{\tau} := \innerprod{\gvsig(\bX)_{st}^{\tau},\ell} = \ell_{\mathbf{1}} +
\sum_{\substack{\vtree\in \fset \\1 \, \leq \, |\vtree| \, \leq \, N}}
\ell_{\vtree}
\innerprod{\gvsig(\bX)_{st}^{\tau},\vtree}.
\]
\end{defn}

\begin{example}
Let \(N=2\) and let \(\cA=\{1,\ldots,d\}\) be the alphabet set as defined before. For \(i,j\in \cA\), we define the compatible Volterra trees \(\vtree^1\), \(\vtree^2\), and \(\vtree^3\) by
\[
\vtree^1 := \left(\Forest{Decorated[[i]]},\,\Forest{Decorated[i]}\right),\qquad
\vtree^2 := \left(\Forest{Decorated[[j[i]]]},\,\Forest{Decorated[j[i]]}\right),\qquad
\vtree^3 := \left(\Forest{Decorated[[i][j]]},\,\Forest{Decorated[j[i]]}\right).
\]
We write \(\vtree^1(i)\), \(\vtree^2(i,j)\), and \(\vtree^3(i,j)\) when indicating the decorations explicitly. By construction, \(|\vtree^1|=1\) and \(|\vtree^2|=|\vtree^3|=2\). The level-two generalized Volterra signature model \(\bV_{\ell}^{2}(\bX)_{st}^{\tau}\) with learnable coefficients indexed by these trees, the forests \(\vtree^1(i)\vtree^1(j)\), and the empty forest \(\mathbf{1}\) is
\[
\begin{aligned}
\bV_{\ell}^{2}(\bX)_{st}^{\tau}
={}&\ell_{\mathbf{1}}+\sum_{i\in\cA}\ell_{\vtree^1(i)}\innerprod{\gvsig(\bX)_{st}^{\tau},\vtree^1(i)}\\
&+\sum_{1\leq i\leq j\leq d}\ell_{\vtree^1(i)\vtree^1(j)}\innerprod{\gvsig(\bX)_{st}^{\tau},\vtree^1(i)}\innerprod{\gvsig(\bX)_{st}^{\tau},\vtree^1(j)}\\
&+\sum_{i,j\in\cA}\left(\ell_{\vtree^2(i,j)}\innerprod{\gvsig(\bX)_{st}^{\tau},\vtree^2(i,j)}+\ell_{\vtree^3(i,j)}\innerprod{\gvsig(\bX)_{st}^{\tau},\vtree^3(i,j)}\right).
\end{aligned}
\]
Equivalently, for a continuous path \(\bX\) of bounded variation and a smooth or at least \(C^1\) kernel \(k\) for which the following integrals are well defined, we have
\[
\begin{aligned}
\bV_{\ell}^{2}(\bX)_{st}^{\tau}
={}&\ell_{\mathbf{1}}+\sum_{i\in\cA}\ell_{\vtree^1(i)}\int_s^t k(\tau,r)\,d\bX_r^{\Forest{Decorated[i]}}\\
&+\sum_{1\leq i\leq j\leq d}\ell_{\vtree^1(i)\vtree^1(j)}\left(\int_s^t k(\tau,r)\,d\bX_r^{\Forest{Decorated[i]}}\right)\left(\int_s^t k(\tau,u)\,d\bX_u^{\Forest{Decorated[j]}}\right)\\
&+\sum_{i,j\in\cA}\ell_{\vtree^2(i,j)}\int_s^t k(\tau,r)\left(\int_s^r k(r,u)\,d\bX_u^{\Forest{Decorated[i]}}\right)d\bX_r^{\Forest{Decorated[j]}}\\
&+\sum_{i,j\in\cA}\ell_{\vtree^3(i,j)}\int_s^t k(\tau,r)\left(\int_s^r k(\tau,u)\,d\bX_u^{\Forest{Decorated[i]}}\right)d\bX_r^{\Forest{Decorated[j]}}.
\end{aligned}
\]
\end{example}

From a computational perspective, we need to establish Chen's identity for the generalized Volterra signature. The main problem lies in the arguments of the product of kernels, where the integration variables are arranged in a chain. To deal with this issue, we need to define generalized convolution product as well.

Let us state generalized Chen's identity and generalized convolution product for the generalized Volterra signature defined in \ref{gen_volterra_def}.

Let \(c\) be an admissible cut of the differential tree \(\dtree\). The corresponding compatible Volterra tree is \(\vtree = (\ktree,\dtree)\). Let \(P^c(\dtree)\) and \(R^c(\dtree)\) be the differential pruned forest and differential trunk respectively. Similarly, let \(P^c(\ktree)\) and \(R^c(\ktree)\) be the kernel pruned forest and kernel trunk respectively. Let the kernel pruned forest be \(P^c(\ktree) = \ktree^1 \dots \ktree^m\) with \(m\) trees. Each tree \(\ktree^j\) for \(j=1,\dots,m\) has a dummy root with no decoration. Let \((P^c(\ktree))^* = (\ktree^1)^* \dots (\ktree^m)^*\) be the forest without the dummy roots and \(\mathbf{r}_1 \dots \mathbf{r}_m\) be the roots of each tree in \((P^c(\ktree))^*\). Set the set \(\bR\) to be the set of these roots i.e., \(\bR = \{\mathbf{r}_1, \dots ,\mathbf{r}_m\}\).

We define the Volterra signature component of a path \(\bX\) corresponding to pruned forest \(P^c(\vtree)\) to be the quantity
\begin{equation} \label{convolution_part}
    \innerprod{\gvsig(\bX)^{\lambda}_{su}, P^c(\vtree)} := \int_{\bA^{P^c(\dtree)}_{su}} \left(\prod_{v\in N_{(P^c(\ktree))^*}} k(r_{\kparent{v}},r_v)\right) \prod_{v\in N_{P^c(\dtree)}} d\bX^v_{r_v},
\end{equation}
where \(\bA^{P^c(\dtree)}_{su}\) is the domain of integration corresponding to differential pruned forest \(P^c(\dtree)\), \(\lambda = (\lambda_1,\dots, \lambda_m)\) are variables corresponding to dummy roots of kernel pruned forest \(P^c(\ktree)\) and 
\[
r_{\kparent{v}} = 
\begin{cases}
    r_{\kparent{v}}, & \text{if } \kparent{v} \in N_{(P^c(\ktree))^*},\\
    \lambda_j \text{ for some } \lambda_j \in \lambda, & \text{otherwise}.
\end{cases}
\]
The signature component \ref{convolution_part} corresponding to \(P^c(\vtree)\) is considered to be a function of free variables \(\lambda = (\lambda_1,\dots, \lambda_m)\). While defining the convolution product, these are replaced by actual integration variables coming from the leaves of the kernel trunk \(R^c(\ktree)\). In other words, \(\lambda_j\) is replaced by the \(r_{\kparent{\mathbf{r}_j}}\), for \(j=1,2,\dots,m\). We define the identification map \(\operatorname{iden}\) that does this replacement as 
\[
\operatorname{iden} (\lambda_j) := r_{\kparent{\mathbf{r}_j}}.
\]
Now, we can define the convolution product as follows.

\begin{defn}[Generalized convolution product] \label{convolution_product}
The convolution product associated to an admissible cut \(c\) is defined as
\begin{align}
    & \innerprod{\gvsig(\bX)^{\tau}_{ut}, 
    R^c(\vtree)} \star \innerprod{\gvsig(\bX)^{\cdot}_{su}, P^c(\vtree)} \\
    & := \int_{\bA^{R^c(\dtree)}_{ut}} \left(\prod_{v\in N_{R^c(\ktree)}} k(r_{\kparent{v}},r_v)\right) \innerprod{\gvsig(\bX)^{\lambda}_{su}, P^c(\vtree)}\Big|_{\operatorname{iden}(\lambda)} \prod_{v\in N_{R^c(\dtree)}} d\bX^v_{r_v},
\end{align}
where \(\bA^{R^c(\dtree)}_{ut}\) is the domain of integration corresponding to differential trunk \(R^c(\dtree)\) and \(\operatorname{iden}(\lambda) := (\operatorname{iden}(\lambda_1),\dots, \operatorname{iden}(\lambda_m))\).
\end{defn}

\begin{example}
    Consider the compatible Volterra tree \(\vtree = \left(\Forest{Decorated[[k[i][j]]]}, \Forest{Decorated[k[i][j]]}\right)\) and the cut \(c = \{(i,k), (j,k)\}\). Corresponding pruned forest and trunk are 
    \[
    P^c(\vtree) = \left(\Forest{Decorated[[i]]} \Forest{Decorated[[j]]}, \Forest{Decorated[i]} \Forest{Decorated[j]}\right), \quad \text{ and } \quad R^c(\vtree) = \left(\Forest{Decorated[[k]]}, \Forest{Decorated[k]}\right)
    \]
    respectively. The signature component corresponding to pruned forest \(P^c(\vtree)\) is
    \[
    \innerprod{\gvsig(\bX)^{\lambda}_{su}, P^c(\vtree)} = \left(\int_s^{u}k(\lambda_1, r_i) d\bX^{\Forest{Decorated[i]}}_{r_i}\right) \left(\int_s^{u}k(\lambda_2, r_j) d\bX^{\Forest{Decorated[j]}}_{r_j}\right),
    \]
    where \(\lambda = (\lambda_1, \lambda_2)\) is the vector of free variables. The identification map gives \(\operatorname{iden}(\lambda_1) = \operatorname{iden}(\lambda_2) = r_k\). Therefore, we have
    \begin{align*}
        & \innerprod{\gvsig(\bX)^{\tau}_{ut}, R^c(\vtree)} \star \innerprod{\gvsig(\bX)^{\lambda}_{su}, P^c(\vtree)} \\
        & \qquad \qquad = \int_u^{t}k(\tau, r_k)  \left(\int_s^{u}k(\lambda_1, r_i) d\bX^{\Forest{Decorated[i]}}_{r_i}\right) \left(\int_s^{u}k(\lambda_2, r_j) d\bX^{\Forest{Decorated[j]}}_{r_j}\right)\Big|_{\operatorname{iden}(\lambda)} d\bX^{\Forest{Decorated[k]}}_{r_k}\\
        & \qquad \qquad = \int_u^{t}k(\tau, r_k)  \left(\int_s^{u}k(r_k, r_i) d\bX^{\Forest{Decorated[i]}}_{r_i}\right) \left(\int_s^{u}k(r_k, r_j) d\bX^{\Forest{Decorated[j]}}_{r_j}\right) d\bX^{\Forest{Decorated[k]}}_{r_k}.
    \end{align*}
\end{example}

\begin{defn}[Generalized Chen's identity]\label{Chen}
Let \(\bX:[0,T] \to \RR^d\) be a path of bounded variation and \(\vtree\) be a compatible Volterra forest. The generalized Volterra signature of \(\bX\) satisfies
\begin{align*}
    \innerprod{\gvsig(\bX)^{\tau}_{st}, \vtree} &  = \innerprod{\gvsig(\bX)^{\tau}_{su}, \vtree} + \innerprod{\gvsig(\bX)^{\tau}_{ut}, \vtree} +\\
    &\sum_{c \, \in \operatorname{Adm}(\dtree)} \innerprod{\gvsig(\bX)^{\tau}_{ut}, R^c(\vtree)} \star \innerprod{\gvsig(\bX)^{\lambda}_{su}, P^c(\vtree)},
\end{align*}
called the generalized Chen's identity.
\end{defn}

\section{Weakly Geometric Volterra Signature}
We can define weakly geometric Volterra signatures by introducing a suitable product on compatible Volterra trees. In the case of weakly geometric signatures, the product of signature components, or equivalently iterated integrals, can be expressed as a sum of signature components/iterated integrals. For this reason, instead of working with forests, we restrict ourselves to trees and, in particular, to ladder trees, since the branched trees can again be represented as sum of ladder trees. This idea will be made precise once we introduce a suitable notion of product.

To begin with, let us restrict ourselves to compatible Volterra trees whose differential tree is a ladder tree. We define \(\cL_{\nu}\) as the linear span of such trees, i.e.,
\[
\cL_{\nu} := \myspan_{\RR}\left\{\vtree=(\ktree,\dtree)\in\tset : \dtree \text{ is a ladder tree} \right\}.
\]

Similar to \(\hopf_{\nu}\), the space \(\cL_{\nu}\) is also graded in the following sense:
\[
\cL_{\nu} = \bigoplus_{n \ge 0} \cL_{\nu}^n,
\]
where
\[
\cL_{\nu}^n := \myspan_{\RR}\{\vtree = (\ktree,\dtree) \in \tset: \deg(\vtree) = n,\ \dtree \text{ is a ladder tree}\}, \quad n \ge 1, \quad \text{and} \quad \cL_{\nu}^0 := \RR.\mathbf{1}.
\]
Moreover, the space \(\cL_{\nu}\) is equipped with the shuffle product defined below.

\begin{defn}[Shuffle product of Volterra trees]\label{Shuffle_prodcut}
For two Volterra trees \(\vtree^1 = (\ktree^1, \dtree^1)\) and \(\vtree^2 = (\ktree^2, \dtree^2)\) in \(\cL_{\nu}\), we define their shuffle product \(\vprod\) by
\[
\vtree^1 \vprod \vtree^2 := \left(\ktree^1 \kprod \ktree^2, \dtree^1 \dprod \dtree^2\right),
\]
where \(\ktree^1 \kprod \ktree^2\) denotes the kernel tree product defined by
\[
\ktree^1 \kprod \ktree^2 := B_{+}\left(B_{-}(\ktree^1)B_{-}(\ktree^2)\right),
\]
and \(\dtree^1 \dprod \dtree^2\) denotes the differential tree product, defined recursively as follows:
\[
\dtree \dprod \mathbf{1} = \mathbf{1} \dprod \dtree = \dtree \quad \text{and} \quad \dtree^1 \dprod \dtree^2 := B^{r_1}_{+}\left(B^{r_1}_{-}(\dtree^1) \dprod \dtree^2\right) + B^{r_2}_{+}\left(\dtree^1 \dprod B^{r_2}_{-}(\dtree^2)\right),
\]
where \(r_1\) and \(r_2\) are the roots of the trees \(\dtree^1\) and \(\dtree^2\), respectively. Here, \(B^{r_k}_{+}(\dtree)\) and \(B^{r_k}_{-}(\dtree)\) denote, respectively, grafting the fixed root \(r_k\) onto the tree \(\dtree\) and removing the fixed root \(r_k\) from \(\dtree\). The product \(\vprod\) is extended bilinearly to all elements of \(\cL_{\nu}\).

By the linear structure of \(\cL_{\nu}\), the shuffle product can equivalently be written as
\[
\vtree^1 \vprod \vtree^2 := \left(\ktree^1 \kprod \ktree^2, B^{r_1}_{+}\left(B^{r_1}_{-}(\dtree^1) \dprod \dtree^2\right)\right) + \left(\ktree^1 \kprod \ktree^2, B^{r_2}_{+}\left(\dtree^1 \dprod B^{r_2}_{-}(\dtree^2)\right)\right).
\]
\end{defn}

\begin{example}
For \(\vtree^1 = \left(\Forest{Decorated[[i]]}, \Forest{Decorated[i]}\right)\) and \(\vtree^2 = \left(\Forest{Decorated[[j]]}, \Forest{Decorated[j]}\right)\), we have
\[
\ktree^1 \kprod \ktree^2 = \Forest{Decorated[[i]]} \kprod \Forest{Decorated[[j]]} = \Forest{Decorated[[i][j]]},
\]
and
\[
\dtree^1 \dprod \dtree^2 = \Forest{Decorated[i]} \dprod \Forest{Decorated[j]} = \Forest{Decorated[j[i]]} + \Forest{Decorated[i[j]]}.
\]
Therefore,
\[
\vtree^1 \vprod \vtree^2 = \left(\Forest{Decorated[[i][j]]}, \Forest{Decorated[j[i]]} + \Forest{Decorated[i[j]]}\right).
\]
\end{example}

Next, we show that the product defined in Definition~\ref{Shuffle_prodcut} is well-defined on \(\cL_{\nu}\) by proving that \(\cL_{\nu}\) is closed under the shuffle product.

\begin{prop}[Closure under the shuffle product]
Let \(\vtree^1 = (\ktree^1,\dtree^1)\) and \(\vtree^2 = (\ktree^2,\dtree^2)\) be two compatible Volterra trees in \(\cL_{\nu}\). Then their shuffle product
\[
\vtree^1 \vprod \vtree^2 := \left(\ktree^1 \kprod \ktree^2, \dtree^1 \dprod \dtree^2\right)
\]
is again an element of \(\cL_{\nu}\). In particular,
\[
\vtree^1,\vtree^2 \in \cL_{\nu} \quad \Longrightarrow \quad \vtree^1 \vprod \vtree^2 \in \cL_{\nu}.
\]
\end{prop}

\begin{proof}
To prove the proposition, we first show that the product of differential trees \(\dtree^1 \dprod \dtree^2\) is a ladder differential tree or a linear combination of ladder differential trees. Secondly, we show that the kernel tree product \(\ktree^1 \kprod \ktree^2\) is compatible with all the differential trees appearing in \(\dtree^1 \dprod \dtree^2\).

We use induction to prove the first claim. Consider two compatible Volterra trees \(\vtree^1, \, \vtree^2 \in \cL_{\nu}\) of degrees \(m\) and \(n\), respectively. If \(m=0\) or \(n=0\), the claim holds trivially since
\[
\dtree^1 \dprod \mathbf{1} = \mathbf{1} \dprod \dtree^1 = \dtree^1 \quad \text{and} \quad \mathbf{1} \dprod \dtree^2 = \dtree^2 \dprod \mathbf{1} = \dtree^2,
\]
where \(\dtree^1\) and \(\dtree^2\) are ladder trees. If \(m=n=1\), i.e., \(\dtree^1 = \Forest{Decorated[i]}\) and \(\dtree^2 = \Forest{Decorated[j]}\) for some decorations \(i\) and \(j\) from the alphabet set, then
\[
\dtree^1 \dprod \dtree^2 = \Forest{Decorated[i]} \dprod \Forest{Decorated[j]} = \Forest{Decorated[j[i]]} + \Forest{Decorated[i[j]]},
\]
which is a linear combination of ladder trees.

Now assume that the product \(\dtree^1 \dprod \dtree^2\) is a linear combination of ladder trees whenever the sum of the degrees of \(\dtree^1\) and \(\dtree^2\) is \(m+n\). We prove that the claim also holds when the sum of the degrees is \(m+n+1\). Without loss of generality, consider a degree--\(m\) tree \(\vtree^1\) and a degree--\((n+1)\) tree \(\vtree^2\). Their differential tree product is defined as
\[
\dtree^1 \dprod \dtree^2 := B^{r_1}_{+}\left(B^{r_1}_{-}(\dtree^1) \dprod \dtree^2\right) + B^{r_2}_{+}\left(\dtree^1 \dprod B^{r_2}_{-}(\dtree^2)\right).
\]
The sum of the degrees of the trees appearing in each of the products \(B^{r_1}_{-}(\dtree^1) \dprod \dtree^2\) and \(\dtree^1 \dprod B^{r_2}_{-}(\dtree^2)\) is \(m+n\). Therefore, by the induction hypothesis, both products are linear combinations of ladder trees. Since adding a root to a ladder tree gives another ladder tree, the right-hand side is again a linear combination of ladder trees. This proves our first claim.

Secondly, let us prove the compatibility of the kernel tree product \(\ktree^1 \kprod \ktree^2\) with all the differential tree appearing in \(\dtree^1 \dprod \dtree^2\). To begin with, observe that the kernel tree product is nothing but an identification of the dummy root i.e., the actual order of the vertices in both of the kernel trees remains same in the product. Also, the order of the vertices is preserved in the product of differential trees as in the product we remove the root and add the same root at the same time. Hence, by using the compatibility of the kernel trees \(\ktree^1\) and \(\ktree^2\) respectively with \(\dtree^1\) and \(\dtree^2\), we have the compatibility of the kernel tree product \(\ktree^1 \kprod \ktree^2\) with all the differential tree appearing in \(\dtree^1 \dprod \dtree^2\).
\end{proof}

Next, we define the weakly geometric Volterra signature. Clearly, we cannot claim that this is also geometric Volterra signature as we are not approximating it using the Volterra signature of the smooth paths. But the two notions should coincide at least for the smooth or \(C^1\) signals.
\begin{defn}[Weakly geometric Volterra signature]
Fix \(0<\gamma<\alpha<1\) and set \(\rho:=\alpha-\gamma\). Here \(\alpha\) represents the H\"older regularity of the driving path, while \(\gamma\) is the singularity of the kernel, corresponding to
\[
|K(\tau,r)| \lesssim |\tau-r|^{-\gamma},
\qquad 0\le r<\tau\le T.
\]
Let \(\cL_\nu\) be the vector space generated by compatible Volterra trees whose differential trees are ladders, including the empty tree \(\mathbf 1\). Set
\[
\Delta^3_{T}
:=
\{(s,t,\tau):0\le s\le t\le\tau\le T\}.
\]
A map
\[
\geovsig(\bX)^{\tau}_{st}:
\Delta^3_{T}\longrightarrow
\cL_\nu,
\qquad
(s,t,\tau)\longmapsto\geovsig(\bX)_{st}^{\tau},
\]
is called a \emph{weakly geometric Volterra signature} if it satisfies the following
properties.
\begin{enumerate}
\item \textbf{Analytic regularity.}
For every basis tree $\vtree\in\cL_\nu$ with
$n:=|\vtree|\ge1$, there exists $C_{\vtree}<\infty$ such that
\[
\left|\innerprod{\geovsig(\bX)_{st}^{\tau},\vtree}\right|
\le C_{\vtree}
\left(
|\tau-t|^{-\gamma}|t-s|^{n\rho+\gamma}
\wedge |\tau-s|^{n\rho}
\right),
\qquad 0\le s<t<\tau\le T,
\]
where $\rho:=\alpha-\gamma>0$.
\item \textbf{Volterra shuffle identity.}
For all \(\vtree^1,\vtree^2\in\mathcal L_\nu\),
\[
\innerprod{\geovsig(\bX)_{st}^{\tau},\mathbf 1}=1,
\]
and
\[
\innerprod{\geovsig(\bX)_{st}^{\tau},\vtree^1}
\innerprod{\geovsig(\bX)_{st}^{\tau},\vtree^2} = \innerprod{
\geovsig(\bX)_{st}^{\tau},
\vtree^1\vprod\vtree^2
}.
\]
\item \textbf{Chen's identity.} For a compatible Volterra tree \(\vtree = (\dtree, \dtree) \in \cL_\nu\)
\begin{align*}
    \innerprod{\geovsig(\bX)^{\tau}_{st}, \vtree} &  = \innerprod{\geovsig(\bX)^{\tau}_{su}, \vtree} + \innerprod{\geovsig(\bX)^{\tau}_{ut}, \vtree} +\\
    &\sum_{c \, \in \operatorname{Adm}(\dtree)} \innerprod{\geovsig(\bX)^{\tau}_{ut}, R^c(\vtree)} \star \innerprod{\geovsig(\bX)^{\lambda}_{su}, P^c(\vtree)},
\end{align*}
for all nontrivial admissible cuts \(c\) of the corresponding differential tree \(\dtree\).
\end{enumerate}
\end{defn}

Let us state and prove the Volterra shuffle identity for the geometric Volterra signature of the path of bounded variation \(\bX:[0,T] \to \RR^d\).

\begin{thm}[Shuffle identity for Volterra signatures]\label{Shuffle_identity}
Let \(\bX:[0,T]\to\mathbb{R}^d\) be a continuous path of bounded variation, and let \(\vtree^1=(\ktree^1,\dtree^1)\) and \(\vtree^2=(\ktree^2,\dtree^2)\) be compatible Volterra trees in \(\cL_{\nu}\). Fix \(0\leq s\leq t\leq\tau\leq T\), and let \(k\) be a bounded measurable real-valued Volterra kernel. For \(i=1,2\), we use \(r_i\) to denote both the root of \(\dtree^i\) and its associated integration variable. Then the paths
\[
G^i_{r_i}:=\innerprod{\geovsig(\bX)_{s r_i}^{\tau},\vtree^i},\qquad s\leq r_i\leq t,\quad i=1,2,
\]
are continuous and of bounded variation, and the following shuffle identity holds:
\[
\innerprod{\geovsig(\bX)_{st}^{\tau},\vtree^1}\innerprod{\geovsig(\bX)_{st}^{\tau},\vtree^2}
=\innerprod{\geovsig(\bX)_{st}^{\tau},\vtree^1\vprod\vtree^2}.
\]
\end{thm}

\begin{proof}
Fix \(0\leq s<t\leq\tau\leq T\), and let \(\vtree^1=(\ktree^1,\dtree^1)\) and \(\vtree^2=(\ktree^2,\dtree^2)\) be two nonempty compatible Volterra trees in \(\cL_{\nu}\). For \(i=1,2\), let \(r_i\) denote the root of \(\dtree^i\), and use \(r_i\) also for its integration time. Recall that \(A_{st}^{\dtree^i}\) is the simplex of integration times ordered by the ladder \(\dtree^i\). We write \(N_{(\ktree^i)^*}\) for the vertices of \(\ktree^i\) other than its external root. In \(k(p_{\ktree^i}(v),v)\), the two arguments are the times of the parent of \(v\) and of \(v\), respectively, with time \(\tau\) assigned to the external root. The superscript \(v\) in \(d\bX^v_{r_v}\) denotes the decoration of \(v\).

\noindent
For \(i=1,2\), define
\[
G^i_t:=\innerprod{\geovsig(\bX)_{st}^{\tau},\vtree^i}
=\int_{A_{st}^{\dtree^i}}
  \prod_{v\in N_{(\ktree^i)^*}} k(p_{\ktree^i}(v),v)
  \prod_{v\in N_{\dtree^i}} d\bX^v_{r_v}.
\]
For a continuous path \(\bX\) of bounded variation and a bounded measurable kernel \(k\), the paths \(G^1\) and \(G^2\) are continuous and of bounded variation.

\noindent
Since \(r_1\) and \(r_2\) are the outermost integration times of their respective ladders, we have
\[
\begin{aligned}
dG^1_{r_1}
&=\left(\int_{A_{sr_1}^{B_-^{r_1}(\dtree^1)}}
  \prod_{v\in N_{(\ktree^1)^*}} k(p_{\ktree^1}(v),v)
  \prod_{v\in N_{\dtree^1}\setminus\{r_1\}} d\bX^v_{r_v}\right)d\bX^{r_1}_{r_1},\\
dG^2_{r_2}
&=\left(\int_{A_{sr_2}^{B_-^{r_2}(\dtree^2)}}
  \prod_{v\in N_{(\ktree^2)^*}} k(p_{\ktree^2}(v),v)
  \prod_{v\in N_{\dtree^2}\setminus\{r_2\}} d\bX^v_{r_v}\right)d\bX^{r_2}_{r_2}.
\end{aligned}
\]
Here \(B_-^{r_i}\) acts only on the ladder tree. The root time \(r_i\) is held fixed in the inner integral, so every kernel factor involving \(r_i\) is still evaluated at that time. Since \(G^1_s=G^2_s=0\), integration by parts gives
\[
G^1_tG^2_t=\int_s^t G^2_{r_1}\,dG^1_{r_1}+\int_s^t G^1_{r_2}\,dG^2_{r_2}.
\]

\noindent
For the first integral, substituting \(G^2_{r_1}\) and \(dG^1_{r_1}\) gives
\[
\begin{aligned}
\int_s^t G^2_{r_1}\,dG^1_{r_1}
={}&\int_s^t
\left(\int_{A_{sr_1}^{\dtree^2}}
  \prod_{v\in N_{(\ktree^2)^*}} k(p_{\ktree^2}(v),v)
  \prod_{v\in N_{\dtree^2}} d\bX^v_{r_v}\right)\\
&\quad\times
\left(\int_{A_{sr_1}^{B_-^{r_1}(\dtree^1)}}
  \prod_{v\in N_{(\ktree^1)^*}} k(p_{\ktree^1}(v),v)
  \prod_{v\in N_{\dtree^1}\setminus\{r_1\}} d\bX^v_{r_v}\right)d\bX^{r_1}_{r_1}.
\end{aligned}
\]
Both external roots carry the same time \(\tau\). Let \(\ktree\) be the kernel tree obtained by identifying these two roots and retaining all their branches. Each non-root vertex then has the same parent time as in its original kernel tree. Thus, when we multiply the two inner integrands, their kernel factors satisfy
\[
\left(\prod_{v\in N_{(\ktree^1)^*}} k(p_{\ktree^1}(v),v)\right)
\left(\prod_{v\in N_{(\ktree^2)^*}} k(p_{\ktree^2}(v),v)\right)
=\prod_{v\in N_{\ktree^*}} k(p_{\ktree}(v),v).
\]
By construction, \(\ktree=B_+\bigl(B_-(\ktree^1)B_-(\ktree^2)\bigr)=\ktree^1\kprod\ktree^2\).

\noindent
For fixed \(r_1\), the two inner domains form the product \(A_{sr_1}^{B_-^{r_1}(\dtree^1)}\times A_{sr_1}^{\dtree^2}\). Splitting this product according to the relative order of the integration times gives one simplex for each interleaving of \(B_-^{r_1}(\dtree^1)\) and \(\dtree^2\) that preserves the order within each ladder. Equal times have zero measure since \(\bX\) is continuous. These interleavings give precisely \(B_-^{r_1}(\dtree^1)\dprod\dtree^2\), and the outer integration places \(r_1\) at the root of each resulting ladder. Hence, by linearity,
\[
\int_s^t G^2_{r_1}\,dG^1_{r_1}
=\innerprod{\geovsig(\bX)_{st}^{\tau},
  \left(\ktree,B_+^{r_1}\bigl(B_-^{r_1}(\dtree^1)\dprod\dtree^2\bigr)\right)}.
\]

\noindent
Similarly, the second integral is
\[
\begin{aligned}
\int_s^t G^1_{r_2}\,dG^2_{r_2}
={}&\int_s^t
\left(\int_{A_{sr_2}^{\dtree^1}}
  \prod_{v\in N_{(\ktree^1)^*}} k(p_{\ktree^1}(v),v)
  \prod_{v\in N_{\dtree^1}} d\bX^v_{r_v}\right)\\
&\quad\times
\left(\int_{A_{sr_2}^{B_-^{r_2}(\dtree^2)}}
  \prod_{v\in N_{(\ktree^2)^*}} k(p_{\ktree^2}(v),v)
  \prod_{v\in N_{\dtree^2}\setminus\{r_2\}} d\bX^v_{r_v}\right)d\bX^{r_2}_{r_2}.
\end{aligned}
\]
The product of the kernel factors is the same as in the first integral. Here the inner domains give the interleavings of \(\dtree^1\) and \(B_-^{r_2}(\dtree^2)\), and \(r_2\) is the outermost integration time. Therefore,
\[
\int_s^t G^1_{r_2}\,dG^2_{r_2}
=\innerprod{\geovsig(\bX)_{st}^{\tau},
  \left(\ktree,B_+^{r_2}\bigl(\dtree^1\dprod B_-^{r_2}(\dtree^2)\bigr)\right)}.
\]

\noindent
Adding the two terms and using the recursive definition of \(\dprod\), we obtain
\[
\begin{aligned}
\innerprod{\geovsig(\bX)_{st}^{\tau},\vtree^1}
\innerprod{\geovsig(\bX)_{st}^{\tau},\vtree^2}
&=\innerprod{\geovsig(\bX)_{st}^{\tau},
  \left(\ktree,\dtree^1\dprod\dtree^2\right)}\\
&=\innerprod{\geovsig(\bX)_{st}^{\tau},\vtree^1\vprod\vtree^2},
\end{aligned}
\]
where the last equality follows from Definition~\ref{Shuffle_prodcut}. If either factor is the unit, the identity follows from the unit convention.
\end{proof}


\section{Uniqueness and Universal Approximation}

In this section we will prove injectivity and the universal approximation for generalized and geometric Volterra signatures. For both types of results we assume that the generalized or the geometric Volterra \(p\)-lift and the corresponding extension of a path exists and is continuous in a suitable topology for some \(p \in \NN\) depending on the regularity of the path. We will divide the discussion into two subsections covering each notion of the Volterra signature.

\subsection{Injectivity and universal approximation for the generalized Volterra signature}
Assume that the full generalized Volterra signature of the path is well defined. Our uniqueness argument concerns time-extended paths and does not require additional regularity beyond that needed for the signature to exist. In what follows, we consider \(\alpha\)-H\"older paths with \(\alpha \in (0,1)\). Uniqueness of the branched signature of time-extended paths is established in \cite{ali2025branched}, while a related result for time-extended semimartingales is proved in \cite{cuchiero2023signature}. Our result extends these uniqueness statements to the generalized Volterra setting, with the corresponding non-Volterra cases recovered by taking the unit kernel, \(k \equiv 1\).

\begin{lem}[Injectivity of the time-augmented generalized Volterra signature] \label{geovsig_uniqueness}
Let \(\bX,\bY:[0,T]\to\mathbb R^d\) be continuous \(\alpha\)-H\"older paths with \(\bX_0=\bY_0=0\) for some \(\alpha \in  (0,1)\).
Construct the time-augmented paths \(\widehat{\bX}(t):=(t,\bX_t)\) and \(\widehat{\bY}(t):=(t,\bY_t)\) in \(\mathbb R^{d+1}\). Assume that the common scalar Volterra kernel \(k\)
satisfies
\[
k(\tau,\cdot)\in C^1([0,T]), \qquad  \tau \ge T, 
\qquad
k(\tau,r)>0
\quad\text{for every }r\in[0,T].
\]
Also, assume that the full generalized Volterra signatures of \(\widehat{\bX}\) and \(\widehat{\bY}\) are well defined and their terminal generalized Volterra signatures coincide at all levels, i.e.,
\[
\gvsig(\widehat{\bX})_{0,T}^{\tau}=\gvsig(\widehat{\bY})_{0,T}^{\tau}.
\]
Then \(\bX_t=\bY_t\) for all \(t\in[0,T]\).
\end{lem}

\begin{proof}
To begin with, we fix a spatial index $i\in \{1,2,\dots,d\}$, set $\bZ^{\Forest{Decorated[i]}}:=\bX^{\Forest{Decorated[i]}}-\bY^{\Forest{Decorated[i]}}$, and define \(\widehat{\bZ} = (t, \bZ)\), a continuous $\alpha$-H\"older path with $\widehat{\bZ}_0=0$.
We use the family of signature coordinates of the path that contain exactly one spatial letter \(i\).
Let the kernel and differential trees be
\[
\ktree = B_+(\underbrace{\Forest{Decorated[0]} \dots \Forest{Decorated[0]} }_{m \text{-times}} \Forest{Decorated[i]} \underbrace{\Forest{Decorated[0]} \dots \Forest{Decorated[0]} }_{k \text{-times}}) \quad \text{ and }  \quad 
\dtree = \underbrace{B_+^{\Forest{Decorated
[0]}} ( \dots( B_+^{\Forest{Decorated
[0]}} }_{k \text{-times}} (B_+^{\Forest{Decorated
[i]}} (\dtree^0) )))),
\]
where \( \dtree^0 = \underbrace{B_+^{\Forest{Decorated
[0]}} ( \dots( B_+^{\Forest{Decorated
[0]}} }_{(m-1)\text{-times}} (\Forest{Decorated[0]})))\). 
We define \(\vtree = (\ktree, \dtree)\). Clearly, the kernel tree without the dummy root is just a forest product of trees with single vertex and the differential tree is just a ladder tree. The component of the signature of \(\widehat{\bZ}\) corresponding to this tree is
\begin{equation}
    \innerprod{\gvsig(\widehat{\bZ})_{0T}^{\tau}, \vtree} = \int_{0 < r_0 < r_1 < \dots < r_k< T} \prod_{i=0}^{k} k(\tau, r_i) \innerprod{\gvsig(\widehat{\bZ})_{0r_0}^{\tau}, \vtree^0} d \widehat{\bZ}_{r_0}^{\Forest{Decorated[i]}} dr_1 \dots dr_k.
\end{equation}
We claim that the simple integration tools lead us to the following representation of the same signature component
\begin{equation} \label{IBP etc.}
    \innerprod{\gvsig(\widehat{\bZ})_{0T}^{\tau}, \vtree} = \frac{1}{m! \, k!} \int_0^T \left(\chi^{\tau}(v)\right)^{m} \left(\chi^{\tau}(T) - \chi^{\tau}(v)\right)^{k} a(v) d \widehat{\bZ}_{v}^{\Forest{Decorated[i]}},
\end{equation}
where \(\chi^{\tau}(v) = \int_0^v a(r)dr\) with \(a(r) = k(\tau,r)\). Let us prove this claim first. First, consider
\begin{equation}
    \Upsilon^{\tau}(u) := \innerprod{\gvsig(\widehat{\bZ})_{0u}^{\tau}, \vtree^0} =  \int_{0 < r_1 < \dots < r_m< u} \prod_{i=1}^{m} k(\tau, r_i) dr_1 \dots dr_m.
\end{equation}
We define \(\Upsilon^{\tau}(u) = 1\) for \(m=0\) and claim that \(\Upsilon^{\tau}(u) = \frac{\left(\chi^{\tau}(u)\right)^{m}}{m!}\). This can be proved using induction easily.
For \(m=1\) the conclusion is obvious. say it is true for \(m=\ell\). Then for \(m=\ell+1\) 
\begin{align*}
    \Upsilon^{\tau}(u) & =  \int_{0}^{u} k(\tau, r_{\ell+1}) \left( \int_{0 < r_1 < \dots < r_{\ell}< r_{\ell+1}} \prod_{i=1}^{\ell} k(\tau, r_i) dr_1 \dots dr_{\ell} \right) dr_{\ell+1}\\
    & =  \int_{0}^{u} k(\tau, r_{\ell+1}) \frac{\left(\chi^{\tau}(r_{\ell+1})\right)^{\ell}}{\ell !} dr_{\ell+1}\\
    & = \frac{1}{\ell!} \int_{0}^{u} \left(\chi^{\tau}(r_{\ell+1})\right)^{\ell} \left(\chi^{\tau}(r_{\ell+1})\right)' dr_{\ell+1}\\
    & = \frac{\left(\chi^{\tau}(u)\right)^{\ell+1}}{(\ell+1)!}
\end{align*}
Let us define \(\Xi^{\tau}_k(t)\) recursively as follows:
\begin{equation}
   \Xi^{\tau}_0(t) =  \int_0^t \Upsilon^{\tau}(u) \, a(u) \, d \widehat{\bZ}_u^{\Forest{Decorated[i]}} \quad \text{ and } \quad \Xi^{\tau}_k(t) =  \int_0^t  \Xi^{\tau}_{k-1}(v)  \dd \chi^{\tau}(v),
\end{equation}
where \(\Xi^{\tau}_0(t)\) is defined in Young's sense because \(\Upsilon^{\tau}(u) \, a(u)\) is \(C^1\). Also, \(\Xi^{\tau}_{\ell}(0) = 0\) for every \(\ell \ge 0\) and \(\Xi^{\tau}_k(T) = \innerprod{\gvsig(\widehat{\bZ})_{0T}^{\tau}, \vtree}\). We define
\begin{equation}
    R_0(v) = 1,  \qquad \text{ and } \qquad  R_j(v) = \frac{\left(\chi^{\tau}(T) - \chi^{\tau}(v)\right)^j}{j!}, \quad j = 1,2, \dots, k,
\end{equation}
which gives
\[
R_{j}(T) = 0,  \text{ for } 1 \le j \le k\quad \text{ and } \quad d R_{j+1}(v) = - R_j(v) d \chi^{\tau}(v), \text{ for } 0 \le j \le k\quad
\]
A simple integration by parts leads to the following relation.
\begin{align} \label{IBP_relation}
    \int_0^T R_j(v) \dd \Xi^{\tau}_k(v) & = \int_0^T R_j(v) \Xi^{\tau}_{k-1}(v) \dd \chi^{\tau}(v) \nonumber \\
    & = -\int_0^T \Xi^{\tau}_{k-1}(v)  \dd R_{j+1}(v) \nonumber \\
    & = \int_0^T R_{j+1}(v) \dd  \Xi^{\tau}_{k-1}(v).
\end{align}
Using relation ~\eqref{IBP_relation} repeatedly, we get the following.
\begin{align} \label{recurrence relation}
    \Xi^{\tau}_k(T) & = \int_0^T R_{1}(v) \dd  \Xi^{\tau}_{k-1}(v) \nonumber \\
    & = \int_0^T R_{2}(v) \dd  \Xi^{\tau}_{k-2}(v) \nonumber \\
    & \qquad \vdots \nonumber\\
    & = \int_0^T R_{k}(v) \dd  \Xi^{\tau}_{0}(v).
\end{align}
Finally, using relation given in ~\eqref{recurrence relation}, the definition of \(\Xi^{\tau}_{0}\), and the Young's substitution rule \cite{castrequini2014youngsystems},  we recover ~\eqref{IBP etc.}.

Moreover, to prove the main claim of the theorem, we use \(k=0\) and arbitrary \(m \ge 0\) in Eq.~\eqref{IBP etc.} to get 
\begin{equation}
    \innerprod{\gvsig(\widehat{\bZ})_{0T}^{\tau}, \vtree} = \frac{1}{m!} \int_0^T \left(\chi^{\tau}(v)\right)^{m} a(v) d \widehat{\bZ}_{v}^{\Forest{Decorated[i]}}.
\end{equation}
By the hypothesis, \(\innerprod{\gvsig(\widehat{\bX})_{0T}^{\tau}, \vtree} = \innerprod{\gvsig(\widehat{\bY})_{0T}^{\tau}, \vtree} \) for every \(\vtree \in \cL_{\nu}\). Therefore
\begin{equation} \label{hypothesis}
    \innerprod{\gvsig(\widehat{\bZ})_{0T}^{\tau}, \vtree} = \frac{1}{m!} \int_0^T \left(\chi^{\tau}(v)\right)^{m} a(v) d \widehat{\bZ}_{v}^{\Forest{Decorated[i]}} = 0.
\end{equation}
If we define \(\bG(t) = \int_0^t a(v) d \widehat{\bZ}_v^{\Forest{Decorated[i]}}\) then ~\eqref{hypothesis} becomes 
\begin{equation*}
 \frac{1}{m!} \int_0^T \left(\chi^{\tau}(v)\right)^{m} d \bG(v) = 0.
\end{equation*}
From the equality of signatures of \(\widehat{\bX}\) and \(\widehat{\bY}\) we can easily establish \(\bG(T) = 0\). Now, a simple Young's integration by parts gives
\begin{equation*}
 -m \int_0^T \bG(v)  \left(\chi^{\tau}(v)\right)^{m-1} a(v) dv = 0, \quad m \ge 1.
\end{equation*}
Therefore,
\begin{equation} \label{main equation}
 \int_0^T \bG(v)  \left(\chi^{\tau}(v)\right)^{n} a(v) dv = 0, \quad n \ge 0.
\end{equation}
Since \(a(v) = k(\tau,v) > 0\), \(\chi^{\tau}(v)\) is strictly increasing, a bijection, and is \(C^1([0,T])\) i.e.,
\begin{equation*}
    \chi^{\tau} : [0,T] \to [0,T'],
\end{equation*}
where \(T' = \chi^{\tau}(T)\). Making the change of variables \(u = \chi^{\tau}(v)\) such that \(du = a(v) dv\) and defining \(\bF(u) : = \bG \left((\chi^{\tau})^{-1}(u)\right)\) in ~\eqref{main equation} gives
\begin{equation*}
    \int_0^{T'} \bF(u) \, u^n du = 0, \quad n \ge 0.
\end{equation*}
By linearity
\begin{equation*}
    \int_0^{T'} \bF(u) \, \mathbf{p}(u) du = 0, \quad \text{ for every polynomial } \mathbf{p}.
\end{equation*}
By invoking the continuity of \(\bF\) and Weierstrass approximation theorem, there exists sequence of polynomials \(\mathbf{p}_j\) such that 
\begin{equation*}
    \|\bF - \mathbf{p}_j\|_{\infty} \to 0 \quad \text{ on } [0,T']. 
\end{equation*}
As a result of this
\begin{equation*}
    \int_0^{T'} \bF(u)^2  \,  du = \lim_{j \to \infty} \int_0^{T'} \bF(u)\, \mathbf{p}_j(u) \, du = 0.
\end{equation*}
The above result together with continuity of \(\bF\) implies \(\bF \equiv 0\), hence \(\bG(t) = 0\), for all \(t \in [0,T]\). Since \(a(v)  = k(\tau,v) > 0\) is \(C^1\) on the compact interval \([0,T]\), so is \(\frac{1}{a(v)}\), therefore
\begin{equation*}
    \widehat{\bZ}_t^{\Forest{Decorated[i]}}  -  \widehat{\bZ}_0^{\Forest{Decorated[i]}}  = \int_0^t \frac{1}{a(v)} \,  a(v) \, d \widehat{\bZ}_v^{\Forest{Decorated[i]}} = \int_0^t \frac{1}{a(v)} \,  d \, \bG(v) = 0
\end{equation*}
Thus \( \widehat{\bZ}^{\Forest{Decorated[i]}} \equiv 0\). Repeating this argument for each \( i \in \{1,2,\dots,d\}\) shows that every spatial component of  \(\widehat{\bZ}\) is zero. Consequently, \(\bZ \equiv 0\) and hence the \(\bX_t = \bY_t\) for all \(t\in [0,T]\).
\end{proof}

Next, we state and prove universal approximation theorem for generalized Volterra signatures.

\begin{thm}[Universal approximation for generalized Volterra signatures of time-extended \(\alpha\)-H\"older paths]
\label{thm:UAT-generalized-Volterra}
Let \(\alpha\in(0,1)\), fix \(\tau\ge T\), and suppose that the kernel \(k\) satisfies the assumptions of Lemma ~\ref{geovsig_uniqueness} for this \(\tau\).

Let \(\hopf_\nu\) be the Hopf algebra of compatible Volterra forests, decorated over the alphabet \(\cS=\{0,1,\dots,d\}\), where \(0\) is the time letter. Denote the canonical pairing between generalized Volterra signatures and \(\hopf_\nu\) by \(\langle\cdot,\cdot\rangle\).

Let \(p=\lfloor \frac{1}{\alpha}\rfloor\). For time-extended paths
\(\bX_t=(t,\bX^-_t)\) equipped with prescribed generalized Volterra lifts, define
\[ \cS^{(p)} := \left\{ \left( \gvsig^p(\bX)_{0t} \right)_{t \in [0,T]}: \bX^-\in C^\alpha([0,T];\RR^d), \quad \bX^-_0=0 \right\}.
\]

Let \(\cH\subset\cS^{(p)}\) be compact, and let \(f:\mathcal H\to\RR\) be continuous. Then, for every \(\varepsilon>0\), there exists \(h\in\hopf_\nu\) such that
\[
\sup_{\left(\gvsig^p(\bX)_{0t}\right)_{t \in [0,T]}\in\mathcal H} \left| f\big(\left(\gvsig^p(\bX)_{0t}\right)_{t \in [0,T]}\big) - \left\langle \gvsig(\bX)_{0,T}^{\tau},h \right\rangle \right| <\varepsilon.
\]
\end{thm}

\begin{proof}
Consider the set
\[
\cG
:=
\left\{
\left(\gvsig^p(\bX)_{0t}\right)_{t\in[0,T]}
\longmapsto
\innerprod{\gvsig(\bX)_{0,T}^{\tau},h}
:
h\in\hopf_\nu
\right\} \subset C(\cH;\RR) .
\]
We need to show that \(\cG\) is a unital subalgebra that separates points. To begin with, \( 1 \in \cG\) because for the unit \(\mathbf{1}\in\hopf_\nu\) we get
\[
\innerprod{\gvsig(\bX)_{0,T}^{\tau},\mathbf{1}}=1.
\]
Thus, \(\cG\) contains the constant functions and, in particular, vanishes nowhere. Moreover, for \(h_1,h_2\in\hopf_\nu\), the Volterra shuffle identity gives
\[
\innerprod{\gvsig(\bX)_{0,T}^{\tau},h_1}
\innerprod{\gvsig(\bX)_{0,T}^{\tau},h_2}
=
\innerprod{
\gvsig(\bX)_{0,T}^{\tau},
h_1 \cdot h_2
}.
\]
Since \(h_1 \cdot h_2\in\hopf_\nu\), the product of any two functions in \(\cG\) also belongs to \(\cG\).
Therefore, \(\cG\) is a unital subalgebra of \(C(\cH;\RR)\).

We now show that \(\cG\) separates points.
Let \(\left(\gvsig^p(\bX)_{0t}\right)_{t\in[0,T]},
\left(\gvsig^p(\bY)_{0t}\right)_{t\in[0,T]}
\in\cH\)
satisfy
\[
\left(\gvsig^p(\bX)_{0t}\right)_{t\in[0,T]}
\ne
\left(\gvsig^p(\bY)_{0t}\right)_{t\in[0,T]}.
\]
Suppose, to the contrary, that
\[
\innerprod{\gvsig(\bX)_{0,T}^{\tau},h}
=
\innerprod{\gvsig(\bY)_{0,T}^{\tau},h}
\qquad
\text{for every }h\in\hopf_\nu.
\]
By Lemma~\ref{geovsig_uniqueness}, we obtain
\(\bX_t=\bY_t\) for any \(t \in [0,T]\).
Since the generalized Volterra lift assigned to each underlying
path is fixed, their prescribed lifts coincide. Consequently,
\[
\left(\gvsig^p(\bX)_{0t}\right)_{t\in[0,T]}
=
\left(\gvsig^p(\bY)_{0t}\right)_{t\in[0,T]},
\]
which is a contradiction.
Therefore, \(\cG\) separates points of \(\cH\). Hence, the claim follows by the real Stone--Weierstrass
theorem.
\end{proof}

\subsection{Injectivity and universal approximation for the geometric Volterra signature}
In this subsection we establish the uniqueness and universal approximation for the geometric Volterra signatures. Again, we assume that \(p\)--lifts of the \(\alpha\)--H\"older paths exist for \( p = \lfloor \frac{1}{\alpha} \rfloor\). Also, the extension to the full geometric Volterra signature exists and is unique.

\begin{lem}[Injectivity of the time-augmented geometric Volterra signature] \label{geovsig_uniqueness}
Let \(\bX,\bY:[0,T]\to\mathbb R^d\) be continuous \(\alpha\)-H\"older paths with \(\bX_0=\bY_0=0\) for some \(\alpha \in  (0,1)\).
Construct the time-augmented paths \(\widehat{\bX}(t):=(t,\bX_t)\) and \(\widehat{\bY}(t):=(t,\bY_t)\) in \(\mathbb R^{d+1}\). Assume that the common scalar Volterra kernel \(k\)
satisfies
\[
k(\tau,\cdot)\in C^1([0,T]), \qquad  \tau \ge T, 
\qquad
k(\tau,r)>0
\quad\text{for every }r\in[0,T].
\]
Also, assume that the full geometric Volterra signatures of \(\widehat{\bX}\) and \(\widehat{\bY}\) are well defined and their terminal geometric Volterra signatures coincide at all levels, i.e.,
\[
\geovsig(\widehat{\bX})_{0,T}^{\tau}=\geovsig(\widehat{\bY})_{0,T}^{\tau}.
\]
Then \(\bX_t=\bY_t\) for all \(t\in[0,T]\).
\end{lem}

\begin{proof}
To begin with, we fix a spatial index $i\in \{1,2,\dots,d\}$, set $\bZ^{\Forest{Decorated[i]}}:=\bX^{\Forest{Decorated[i]}}-\bY^{\Forest{Decorated[i]}}$, and define \(\widehat{\bZ} = (t, \bZ)\), a continuous $\alpha$-H\"older path with $\widehat{\bZ}_0=0$.
We use the family of signature coordinates of the path that contain exactly one spatial letter \(i\).
Let the kernel and differential trees be
\[
\ktree = B_+(\underbrace{\Forest{Decorated[0]} \dots \Forest{Decorated[0]} }_{m \text{-times}} \Forest{Decorated[i]} \underbrace{\Forest{Decorated[0]} \dots \Forest{Decorated[0]} }_{k \text{-times}}) \quad \text{ and }  \quad 
\dtree = \underbrace{B_+^{\Forest{Decorated
[0]}} ( \dots( B_+^{\Forest{Decorated
[0]}} }_{k \text{-times}} (B_+^{\Forest{Decorated
[i]}} (\dtree^0) )))),
\]
where \( \dtree^0 = \underbrace{B_+^{\Forest{Decorated
[0]}} ( \dots( B_+^{\Forest{Decorated
[0]}} }_{(m-1)\text{-times}} (\Forest{Decorated[0]})))\). 
We define \(\vtree = (\ktree, \dtree)\). Clearly, the kernel tree without the dummy root is just a forest product of trees with single vertex and the differential tree is just a ladder tree. The component of the signature of \(\widehat{\bZ}\) corresponding to this tree is
\begin{equation}
    \innerprod{\geovsig(\widehat{\bZ})_{0T}^{\tau}, \vtree} = \int_{0 < r_0 < r_1 < \dots < r_k< T} \prod_{i=0}^{k} k(\tau, r_i) \innerprod{\geovsig(\widehat{\bZ})_{0r_0}^{\tau}, \vtree^0} d \widehat{\bZ}_{r_0}^{\Forest{Decorated[i]}} dr_1 \dots dr_k.
\end{equation}
We claim that the simple integration tools lead us to the following representation of the same signature component
\begin{equation} \label{IBP etc.}
    \innerprod{\geovsig(\widehat{\bZ})_{0T}^{\tau}, \vtree} = \frac{1}{m! \, k!} \int_0^T \left(\chi^{\tau}(v)\right)^{m} \left(\chi^{\tau}(T) - \chi^{\tau}(v)\right)^{k} a(v) d \widehat{\bZ}_{v}^{\Forest{Decorated[i]}},
\end{equation}
where \(\chi^{\tau}(v) = \int_0^v a(r)dr\) with \(a(r) = k(\tau,r)\). Let us prove this claim first. First, consider
\begin{equation}
    \Upsilon^{\tau}(u) := \innerprod{\geovsig(\widehat{\bZ})_{0u}^{\tau}, \vtree^0} =  \int_{0 < r_1 < \dots < r_m< u} \prod_{i=1}^{m} k(\tau, r_i) dr_1 \dots dr_m.
\end{equation}
We define \(\Upsilon^{\tau}(u) = 0\) for \(m=0\) and claim that \(\Upsilon^{\tau}(u) = \frac{\left(\chi^{\tau}(u)\right)^{m}}{m!}\). This can be proved using induction easily.
For \(m=1\) the conclusion is obvious. say it is true for \(m=\ell\). Then for \(m=\ell+1\) 
\begin{align*}
    \Upsilon^{\tau}(u) & =  \int_{0}^{u} k(\tau, r_{\ell+1}) \left( \int_{0 < r_1 < \dots < r_{\ell}< r_{\ell+1}} \prod_{i=1}^{\ell} k(\tau, r_i) dr_1 \dots dr_{\ell} \right) dr_{\ell+1}\\
    & =  \int_{0}^{u} k(\tau, r_{\ell+1}) \frac{\left(\chi^{\tau}(r_{\ell+1})\right)^{\ell}}{\ell !} dr_{\ell+1}\\
    & = \frac{1}{\ell!} \int_{0}^{u} \left(\chi^{\tau}(r_{\ell+1})\right)^{\ell} \left(\chi^{\tau}(r_{\ell+1})\right)' dr_{\ell+1}\\
    & = \frac{\left(\chi^{\tau}(u)\right)^{\ell+1}}{(\ell+1)!}
\end{align*}
Let us define \(\Xi^{\tau}_k(t)\) recursively as follows:
\begin{equation}
   \Xi^{\tau}_0(t) =  \int_0^t \Upsilon^{\tau}(u) \, a(u) \, d \widehat{\bZ}_u^{\Forest{Decorated[i]}} \quad \text{ and } \quad \Xi^{\tau}_k(t) =  \int_0^t  \Xi^{\tau}_{k-1}(v)  \dd \chi^{\tau}(v),
\end{equation}
where \(\Xi^{\tau}_0(t)\) is defined in Young's sense because \(\Upsilon^{\tau}(u) \, a(u)\) is \(C^1\). Also, \(\Xi^{\tau}_{\ell}(0) = 0\) for every \(\ell \ge 0\) and \(\Xi^{\tau}_k(T) = \innerprod{\geovsig(\widehat{\bZ})_{0T}^{\tau}, \vtree}\). We define
\begin{equation}
    R_0(v) = 1,  \qquad \text{ and } \qquad  R_j(v) = \frac{\left(\chi^{\tau}(T) - \chi^{\tau}(v)\right)^j}{j!}, \quad j = 1,2, \dots, k,
\end{equation}
which gives
\[
R_{j}(T) = 0,  \text{ for } 1 \le j \le k\quad \text{ and } \quad d R_{j+1}(v) = - R_j(v) d \chi^{\tau}(v), \text{ for } 0 \le j \le k\quad
\]
A simple integration by parts leads to the following relation.
\begin{align} \label{IBP_relation}
    \int_0^T R_j(v) \dd \Xi^{\tau}_k(v) & = \int_0^T R_j(v) \Xi^{\tau}_{k-1}(v) \dd \chi^{\tau}(v) \nonumber \\
    & = -\int_0^T \Xi^{\tau}_{k-1}(v)  \dd R_{j+1}(v) \nonumber \\
    & = \int_0^T R_{j+1}(v) \dd  \Xi^{\tau}_{k-1}(v).
\end{align}
Using relation ~\eqref{IBP_relation} repeatedly, we get the following.
\begin{align} \label{recurrence relation}
    \Xi^{\tau}_k(T) & = \int_0^T R_{1}(v) \dd  \Xi^{\tau}_{k-1}(v) \nonumber \\
    & = \int_0^T R_{2}(v) \dd  \Xi^{\tau}_{k-2}(v) \nonumber \\
    & \qquad \vdots \nonumber\\
    & = \int_0^T R_{k}(v) \dd  \Xi^{\tau}_{0}(v).
\end{align}
Finally, using relation given in ~\eqref{recurrence relation}, the definition of \(\Xi^{\tau}_{0}\), and the Young's substitution rule \cite{castrequini2014youngsystems},  we recover ~\eqref{IBP etc.}.

Moreover, to prove the main claim of the theorem, we use \(k=0\) and arbitrary \(m \ge 0\) in Eq.~\eqref{IBP etc.} to get 
\begin{equation}
    \innerprod{\geovsig(\widehat{\bZ})_{0T}^{\tau}, \vtree} = \frac{1}{m!} \int_0^T \left(\chi^{\tau}(v)\right)^{m} a(v) d \widehat{\bZ}_{v}^{\Forest{Decorated[i]}}.
\end{equation}
By the hypothesis, \(\innerprod{\geovsig(\widehat{\bX})_{0T}^{\tau}, \vtree} = \innerprod{\geovsig(\widehat{\bY})_{0T}^{\tau}, \vtree} \) for every \(\vtree \in \cL_{\nu}\). Therefore
\begin{equation} \label{hypothesis}
    \innerprod{\geovsig(\widehat{\bZ})_{0T}^{\tau}, \vtree} = \frac{1}{m!} \int_0^T \left(\chi^{\tau}(v)\right)^{m} a(v) d \widehat{\bZ}_{v}^{\Forest{Decorated[i]}} = 0.
\end{equation}
If we define \(\bG(t) = \int_0^t a(v) d \widehat{\bZ}_v^{\Forest{Decorated[i]}}\) then ~\eqref{hypothesis} becomes 
\begin{equation*}
 \frac{1}{m!} \int_0^T \left(\chi^{\tau}(v)\right)^{m} d \bG(v) = 0.
\end{equation*}
From the equality of signatures of \(\widehat{\bX}\) and \(\widehat{\bY}\) we can easily establish \(\bG(T) = 0\). Now, a simple Young's integration by parts gives
\begin{equation*}
 -m \int_0^T \bG(v)  \left(\chi^{\tau}(v)\right)^{m-1} a(v) dv = 0, \quad m \ge 1.
\end{equation*}
Therefore,
\begin{equation} \label{main equation}
 \int_0^T \bG(v)  \left(\chi^{\tau}(v)\right)^{n} a(v) dv = 0, \quad n \ge 0.
\end{equation}
Since \(a(v) = k(\tau,v) > 0\), \(\chi^{\tau}(v)\) is strictly increasing, a bijection, and is \(C^1([0,T])\) i.e.,
\begin{equation*}
    \chi^{\tau} : [0,T] \to [0,T'],
\end{equation*}
where \(T' = \chi^{\tau}(T)\). Making the change of variables \(u = \chi^{\tau}(v)\) such that \(du = a(v) dv\) and defining \(\bF(u) : = \bG \left((\chi^{\tau})^{-1}(u)\right)\) in ~\eqref{main equation} gives
\begin{equation*}
    \int_0^{T'} \bF(u) \, u^n du = 0, \quad n \ge 0.
\end{equation*}
By linearity
\begin{equation*}
    \int_0^{T'} \bF(u) \, \mathbf{p}(u) du = 0, \quad \text{ for every polynomial } \mathbf{p}.
\end{equation*}
By invoking the continuity of \(\bF\) and Weierstrass approximation theorem, there exists sequence of polynomials \(\mathbf{p}_j\) such that 
\begin{equation*}
    \|\bF - \mathbf{p}_j\|_{\infty} \to 0 \quad \text{ on } [0,T']. 
\end{equation*}
As a result of this
\begin{equation*}
    \int_0^{T'} \bF(u)^2  \,  du = \lim_{j \to \infty} \int_0^{T'} \bF(u)\, \mathbf{p}_j(u) \, du = 0.
\end{equation*}
The above result together with continuity of \(\bF\) implies \(\bF \equiv 0\), hence \(\bG(t) = 0\), for all \(t \in [0,T]\). Since \(a(v)  = k(\tau,v) > 0\) is \(C^1\) on the compact interval \([0,T]\), so is \(\frac{1}{a(v)}\), therefore
\begin{equation*}
    \widehat{\bZ}_t^{\Forest{Decorated[i]}}  -  \widehat{\bZ}_0^{\Forest{Decorated[i]}}  = \int_0^t \frac{1}{a(v)} \,  a(v) \, d \widehat{\bZ}_v^{\Forest{Decorated[i]}} = \int_0^t \frac{1}{a(v)} \,  d \, \bG(v) = 0
\end{equation*}
Thus \( \widehat{\bZ}^{\Forest{Decorated[i]}} \equiv 0\). Repeating this argument for each \( i \in \{1,2,\dots,d\}\) shows that every spatial component of  \(\widehat{\bZ}\) is zero. Consequently, \(\bZ \equiv 0\) and hence the \(\bX_t = \bY_t\) for all \(t\in [0,T]\).
\end{proof}

Next, we state and prove universal approximation theorem for geometric Volterra signatures.

\begin{thm}[Universal approximation for geometric Volterra signatures of time-extended $\alpha$-H\"older paths]
\label{thm:UAT-geometric-Volterra}
Let \(\alpha\in(0,1)\), fix \(\tau\ge T\), and suppose that the kernel \(k\) satisfies the assumptions of Lemma ~\ref{geovsig_uniqueness} for this \(\tau\).

Let \(\cL_\nu\) be the unital shuffle algebra of compatible Volterra trees whose differential trees are ladders, decorated over the alphabet \(\cS=\{0,1,\dots,d\}\), where \(0\) is the time letter. Denote the canonical pairing between geometric Volterra signatures and \(\cL_\nu\) by \(\langle\cdot,\cdot\rangle\).

Let \(p=\lfloor \frac{1}{\alpha}\rfloor\). For time-extended paths
\(\bX_t=(t,\bX^-_t)\) equipped with prescribed geometric Volterra lifts, define
\[ \cS^{(p)} := \left\{ \left( \geovsig^p(\bX)_{0t} \right)_{t \in [0,T]}: \bX^-\in C^\alpha([0,T];\RR^d), \quad \bX^-_0=0 \right\}.
\]

Let \(\cH\subset\cS^{(p)}\) be compact, and let \(f:\mathcal H\to\RR\) be continuous. Then, for every \(\varepsilon>0\), there exists \(h\in\cL_\nu\) such that
\[
\sup_{\left(\geovsig^p(\bX)_{0t}\right)_{t \in [0,T]}\in\mathcal H} \left| f\big(\left(\geovsig^p(\bX)_{0t}\right)_{t \in [0,T]}\big) - \left\langle \geovsig(\bX)_{0,T}^{\tau},h \right\rangle \right| <\varepsilon.
\]
\end{thm}

\begin{proof}
Consider the set
\[
\cG
:=
\left\{
\left(\geovsig^p(\bX)_{0t}\right)_{t\in[0,T]}
\longmapsto
\innerprod{\geovsig(\bX)_{0,T}^{\tau},h}
:
h\in\cL_\nu
\right\} \subset C(\cH;\RR) .
\]
We need to show that \(\cG\) is a unital subalgebra that separates points. To begin with, \( 1 \in \cG\) because for the unit \(\mathbf{1}\in\cL_\nu\) we get
\[
\innerprod{\geovsig(\bX)_{0,T}^{\tau},\mathbf{1}}=1.
\]
Thus, \(\cG\) contains the constant functions and, in particular, vanishes nowhere. Moreover, for \(h_1,h_2\in\cL_\nu\), the Volterra shuffle identity gives
\[
\innerprod{\geovsig(\bX)_{0,T}^{\tau},h_1}
\innerprod{\geovsig(\bX)_{0,T}^{\tau},h_2}
=
\innerprod{
\geovsig(\bX)_{0,T}^{\tau},
h_1 \vprod h_2
}.
\]
Since \(h_1 \vprod h_2\in\cL_\nu\), the product of any two functions in \(\cG\) also belongs to \(\cG\).
Therefore, \(\cG\) is a unital subalgebra of \(C(\cH;\RR)\).

We now show that \(\cG\) separates points.
Let \(\left(\geovsig^p(\bX)_{0t}\right)_{t\in[0,T]},
\left(\geovsig^p(\bY)_{0t}\right)_{t\in[0,T]}
\in\cH\)
satisfy
\[
\left(\geovsig^p(\bX)_{0t}\right)_{t\in[0,T]}
\ne
\left(\geovsig^p(\bY)_{0t}\right)_{t\in[0,T]}.
\]
Suppose, to the contrary, that
\[
\innerprod{\geovsig(\bX)_{0,T}^{\tau},h}
=
\innerprod{\geovsig(\bY)_{0,T}^{\tau},h}
\qquad
\text{for every }h\in\cL_\nu.
\]
By Lemma~\ref{geovsig_uniqueness}, we obtain
\(\bX_t=\bY_t\) for any \(t \in [0,T]\).
Since the geometric Volterra lift assigned to each underlying
path is fixed, their prescribed lifts coincide. Consequently,
\[
\left(\geovsig^p(\bX)_{0t}\right)_{t\in[0,T]}
=
\left(\geovsig^p(\bY)_{0t}\right)_{t\in[0,T]},
\]
which is a contradiction.
Therefore, \(\cG\) separates points of \(\cH\). Hence, the claim follows by the real Stone--Weierstrass
theorem.
\end{proof}

\section{Special Cases of the Generalized Volterra Signature}
In this section, we discuss the special cases obtained by imposing restrictions on the Volterra tree \(\vtree=(\ktree,\dtree)\) and the kernel \(k\). We describe the corresponding signatures and their underlying algebraic structures. We first consider restrictions on the kernel tree and the differential tree. We then consider the constant-kernel case, where the kernel factors disappear and the construction reduces to signatures indexed by the differential tree alone.

\subsection{Restrictions on the Volterra tree}
We first consider arbitrary compatible Volterra trees \(\vtree=(\ktree,\dtree)\), without any further restriction. This gives the generalized Volterra signature \(\gvsig(\bX)\), whose underlying algebra is the Hopf algebra \(\hopf_{\nu}\) introduced in this work. Both the kernel tree and the differential tree are retained in the description of each coordinate. The kernel tree records the arrangement of the kernel factors, while the differential tree records the integration structure. An example of a compatible Volterra tree in this case is
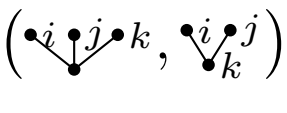
\[\vtree=\left(\Forest{Decorated[[i][j][k]]},\Forest{Decorated[k[i][j]]}\right).\]
The remaining cases in this subsection are obtained by restricting this family of compatible Volterra trees.

If we restrict the differential tree \(\dtree\) to be a ladder tree, we obtain the geometric Volterra signature \(\geovsig(\bX)\). The kernel tree \(\ktree\) can still vary, provided that the pair \((\ktree,\dtree)\) remains compatible. Thus, different compatible kernel trees can occur with the same ladder differential tree. These coordinates correspond to the iterated integrals appearing in \cite{feng2023cubature}. The underlying algebra is \(\cL_{\nu}\), equipped with the product and unit introduced in this work. For example, this case includes
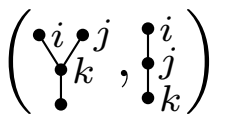
\[\vtree=\left(\Forest{Decorated[[k[i][j]]]},\Forest{Decorated[k[j[i]]]}\right).\]
In this example, the differential tree is a ladder, while the kernel tree has branching. This shows that the ladder restriction is imposed only on the differential tree.

Another restriction is to require \(\ktree=B_+(\dtree)\), so that the kernel tree is obtained by adding a new root above the root of the differential tree. Here the differential tree may still have branching, but once it is chosen, the kernel tree is determined. Apart from the additional root, the two trees therefore have the same branching structure. This gives the iterated integrals associated with the Volterra-type rough paths considered in \cite{bruned2023ramification}. Their algebraic structure is described by the planted rooted-tree Hopf algebra in that work, with its forest product and coproduct. An example is
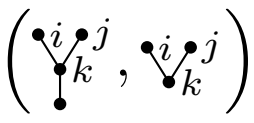
\[\vtree=\left(\Forest{Decorated[[k[i][j]]]},\Forest{Decorated[k[i][j]]}\right).\]

Finally, if \(\dtree\) is a ladder tree and \(\ktree=B_+(\dtree)\), we recover the scalar-kernel case of the classical Volterra signature considered in \cite{hager2026volterra}. Thus, this case imposes both of the preceding restrictions. The differential tree is a ladder, and the kernel tree is the corresponding planted ladder tree. For example,
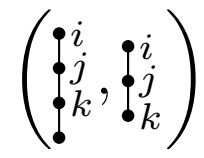
\[\vtree=\left(\Forest{Decorated[[k[j[i]]]]},\Forest{Decorated[k[j[i]]]}\right).\]
This is a more restrictive choice than the geometric Volterra case, where the kernel tree is allowed to have branching. As recalled earlier, the classical Volterra signature satisfies a Chen identity involving the Volterra convolution. However, the usual word-shuffle identity does not hold in general for its coordinates.

\subsection{The constant-kernel case}
We now impose the additional restriction \(k\equiv1\) and take \(\ktree=B_+(\dtree)\). Since every kernel factor equals \(1\), the kernel tree can be omitted, leaving only the differential tree \(\dtree\). The resulting iterated integrals therefore carry no kernel weights. The distinction between the following two cases is whether the differential tree is allowed to have branching or is restricted to a ladder tree.

If \(\dtree\) is an arbitrary rooted tree, the generalized Volterra signature reduces to the branched signature \(\bsig(\bX)\) considered in \cite{gubinelli2010ramification}. Its coordinates are indexed by rooted trees, and their products are represented by forests. Branching in the differential tree records the products that appear inside the recursive definition of the iterated integrals. By retaining all rooted trees, we retain these products of iterated integrals. The underlying algebra is the Butcher--Connes--Kreimer Hopf algebra.

If we further restrict \(\dtree\) to be a ladder tree, we recover the classical signature \(\sig(\bX)\). The vertices of a ladder tree form a single chain, so their decorations can be read in the order of integration to give a word. The corresponding coordinates are therefore the usual word-indexed iterated integrals. Products of these integrals satisfy the classical shuffle identity, and the underlying algebra is the shuffle Hopf algebra; see \cite{cass2022combinatorial}.

The restrictions on the differential tree \(\dtree\), the kernel tree \(\ktree\), and the kernel function \(k\) describe the signature coordinates and their formal algebraic structure. These restrictions do not, by themselves, impose an integration convention. The convention must be specified separately, and the algebraic identities satisfied by the resulting integrals must be checked for that choice.

\begin{center}
\resizebox{0.8\linewidth}{!}{%
\begin{tikzpicture}[
    box/.style={
        draw=diagramblue!65,
        line width=0.65pt,
        rounded corners=2pt,
        text=diagramink,
        align=center,
        text width=5.05cm,
        minimum height=1.85cm,
        inner sep=7pt,
        font=\small
    },
    restriction/.style={
        -{Stealth[length=2.2mm,width=1.7mm]},
        draw=diagramink!85,
        line width=0.7pt
    },
    arrowtext/.style={
        font=\small,
        text=diagramink,
        fill=white,
        inner sep=2.5pt,
        align=center
    }
]

\node[box,fill=diagramblue!5] (general) at (0,0)
    {\textbf{Generalized Volterra}\\
     \textbf{signature}\\[4pt]
     {\footnotesize Current work}};

\node[box,fill=diagramblue!5] (geometric) at (8.8,0)
    {\textbf{Geometric Volterra}\\
     \textbf{signature}\\[4pt]
     {\footnotesize Current work and the iterated integrals appearing in \cite{feng2023cubature}}};

\node[box,fill=black!2] (branchedvolterra) at (0,-3.6)
    {\textbf{Volterra-type rough path}\\
     \textbf{signature}\\[4pt]
     {\footnotesize Volterra-type rough path as defined in \cite{bruned2023ramification}}};

\node[box,fill=black!2] (classicalvolterra) at (8.8,-3.6)
    {\textbf{Classical Volterra}\\
     \textbf{signature}\\[4pt]
     {\footnotesize Scalar kernel version of the Volterra signature \cite{hager2026volterra} }};

\node[box,fill=diagramblue!3] (branched) at (0,-7.2)
    {\textbf{Branched signature}\\[4pt]
     {\footnotesize The branched signature of Gubinelli \cite{gubinelli2004controlling}}};

\node[box,fill=diagramblue!3] (classical) at (8.8,-7.2)
    {\textbf{The Classical signature}\\[4pt]
     {\footnotesize The siganture arising from geometric rough paths  \cite{Friz_Victoir_2010} }};

\draw[restriction] (general.east)
    -- node[arrowtext,above] {\(\dtree\) is a ladder}
    (geometric.west);

\draw[restriction] (branchedvolterra.east)
    -- node[arrowtext,above] {\(\dtree\) is a ladder}
    (classicalvolterra.west);

\draw[restriction] (branched.east)
    -- node[arrowtext,above] {\(\dtree\) is a ladder}
    (classical.west);

\draw[restriction] (general.south)
    -- node[arrowtext,right] {\(\ktree=B_+(\dtree)\)}
    (branchedvolterra.north);

\draw[restriction] (geometric.south)
    -- node[arrowtext,left] {\(\ktree=B_+(\dtree)\)}
    (classicalvolterra.north);

\draw[restriction] (branchedvolterra.south)
    -- node[arrowtext,right] {\(k\equiv1\)}
    (branched.north);

\draw[restriction] (classicalvolterra.south)
    -- node[arrowtext,left] {\(k\equiv1\)}
    (classical.north);

\end{tikzpicture}%
}

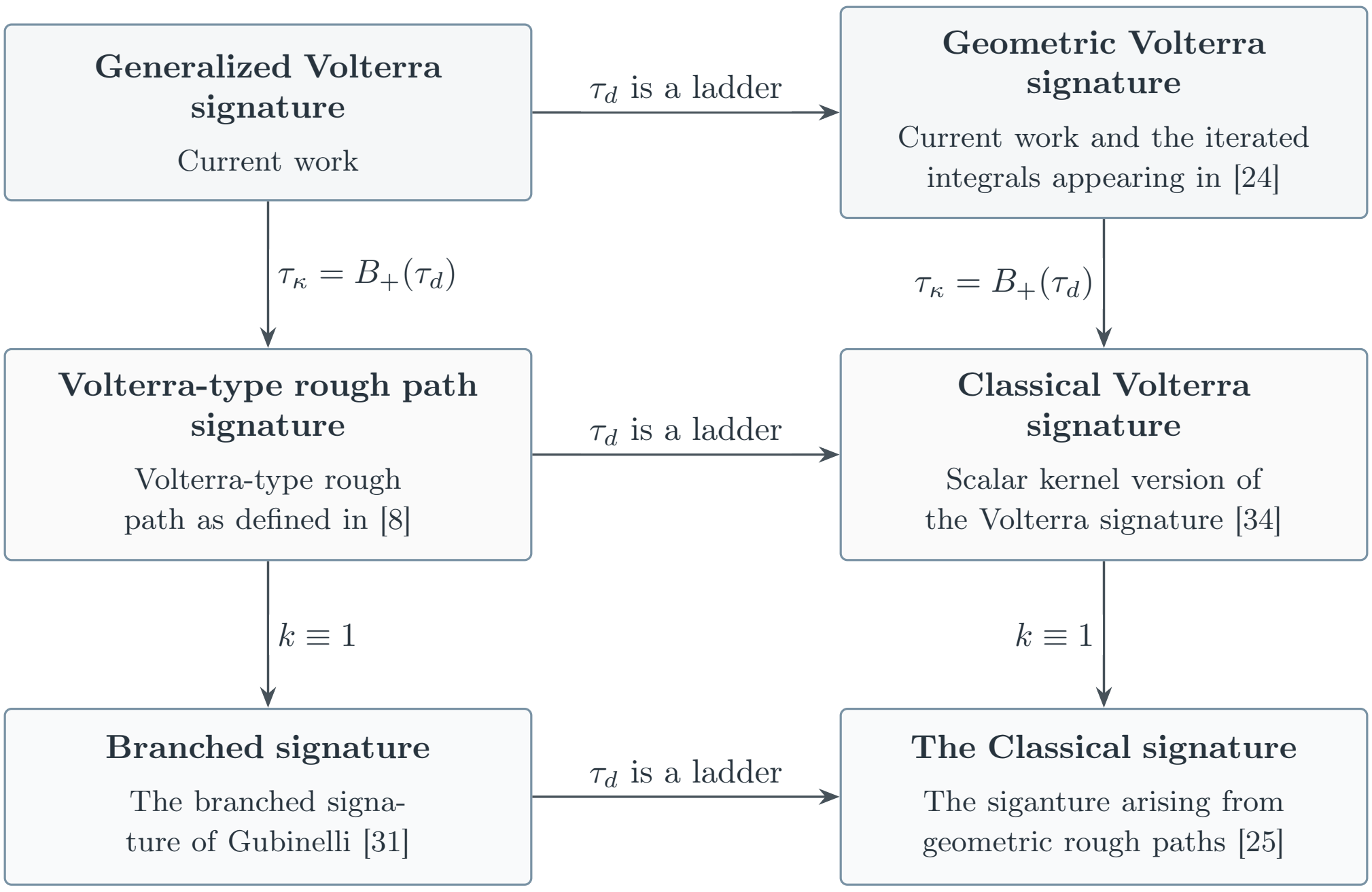
\captionof{figure}{Special cases of the generalized Volterra signature. Each arrow indicates the restriction on the trees or the kernel that gives the corresponding signature.}
\label{fig:volterra-special-cases}
\end{center}

\section{Reduction of the geometric Volterra signature to the classical signature of an extended path}
In this section we provide a method to represent the geometric Volterra signature \(\geovsig\) of a path \(\bX:[0,T] \to \RR^d\) as the classical signature \(\sig\) of an extended path \(\Bar{\bX}\) for special class of kernels. The kernels that we consider in our discussion are separable kernels and are defined below. For the general kernels we will provide an approximation argument later.
\begin{defn}[Separable Kernel] \label{separable kernel}
    A kernel \(k: \Delta^2_T \to \RR\) is said to be separable if it can be decomposed into the form
    \[
    k(t,s) = \frac{g(t)h(s)}{G}, \quad \text{ for } 0 \le s < t \le T, \text{ and } G \ne 0.
    \]
for some smooth functions \(g\) and \(h\). 
\end{defn}
Such type of kernels are considered in \cite{alqahtani2023solution} to study a class of Volterra integral equations. One e example of such kernel  is the exponential kernel, i.e., \(k(t,s) = e^{-\lambda(t-s)}\). Here \(g(t) = e^{-\lambda \, t}\), \(h(s) = e^{\lambda \, s}\), and \(G = 1\). 

For \(0 \le s < t \le \tau \le T\) with \(\tau\) fixed and \(\bX:[0,T]\) be a path of bounded variation. We define the component of the geometric Volterra signature \(\geovsig\) corresponding to the \(\vtree = (\ktree, \dtree)  \in \cL_{\nu}\), where \(\dtree\) is the ladder differential tree say \(\dtree = \Forest{Decorated[i_n[{\scalebox{1}[0.3]{\ensuremath{\vdots}}}[i_2[i_1]]]]}\), as follows
\begin{equation}
\innerprod{\geovsig(\bX)^{\tau}_{st}, \vtree} := \int_{s<r_1< \cdots < r_n < t
} \left(\prod^n_{v=1} k (r_{\kparent{i_v}}, r_{i_v})\right) d\bX_{r_{i_1}}^{i_1} \cdots d\bX_{r_{i_n}}^{i_n},
\end{equation}
where \(k(s,t)\) is the kernel function and \(\kparent{i_v}\) is the parent of vertex \(i_v\) in the kernel tree \(\ktree\).

Before presenting the main results of the section let use rewrite the kernel product \(\prod^n_{v=1} k (r_{\kparent{i_v}}, r_{i_v})\) using the function \(g\) and \(h\) as given in Definition \eqref{separable kernel}. 

\begin{lem}[Kernel product] \label{kernel product}
    Let the kernel \(k(t,s)\) satisfies Definition \eqref{separable kernel}. The kernel product \(\prod^n_{v=1} k (r_{\kparent{i_v}}, r_{i_v})\) can equivalently be written as follows
    \begin{equation} \label{factorization of kernel}
        \prod^n_{v=1} k (r_{\kparent{i_v}}, r_{i_v}) = \frac{g(\tau)^m}{G^n} \prod^n_{v=1} g(r_{i_v})^{c_{i_v}} h(r_{i_v})
    \end{equation}
    where \(m\) is the number of children of the dummy root corresponding to \(\tau\) in the kernel tree \(\ktree\), \(n : = |\vtree|\), and \(c_{i_v}\) is the number of children of vertex \(i_v\) in the kernel tree.
\end{lem}

\begin{proof}
    Consider the kernel tree \(\ktree\) with cardinality of the vertices except the dummy root to be \(n\). Each component in the product \(k (r_{\kparent{i_v}}, r_{i_v})\) can be written using Definition  \eqref{separable kernel} as \( \frac{g (r_{\kparent{i_v}}) h(r_{i_v})}{G}\). There are \(n\) of these contributions. Therefore the terms \(G^n\) and \(h(r_{i_v})\) in ~\eqref{factorization of kernel} make sense. If there are \(m\) immediate children of the dummy root then there are exactly \(m\) vertices whose parent is the dummy root in the kernel tree. Therefore the term \(g (r_{\kparent{i_v}})\) should be \(g(\tau)\), and it happens exactly \(m\) times. This contributes to the term \(g (\tau)^m\). Finally, if \(i_v\) is a vertex with \(c_{i_v}\) many children then for each child the term \(g(r_{i_v})\) appears. This term appears exactly \(c_{i_v}\) many times contributing to the term \(g (r_{{i_v}})^{c_{i_v}}\) in ~\eqref{factorization of kernel}.
\end{proof}
\begin{rem}
    The number of children of the dummy root in the kernel tree and number of children of each vertex \(i_v\) should sum to total number of vertices, i.e.,
    \[
    m + \sum_{v=1}^n c_{i_v} = n.
    \]
\end{rem}
\begin{example}
    Consider the kernel tree \(\ktree =  \Forest{Decorated[[k][j[i]]]}\). There are exactly two children of the dummy root and vertices \(\Forest{Decorated[i]}\), \(\Forest{Decorated[j]}\), and \(\Forest{Decorated[k]}\) have children \(0,1\), and \(0\) respectively. Therefore 
    \begin{align*}
        \frac{g(\tau)^2}{G^3} \prod^3_{v=1} g(r_{i_v})^{c_{i_v}} h(r_{i_v}) & = \frac{g(\tau)^2}{G^3} \left[g(r_i)^0h(r_i) \cdot g(r_j)^1h(r_j) \cdot g(r_k)^0h(r_k)\right]\\
        & = \frac{g(\tau)^2}{G^3} \left[h(r_i) g(r_j)h(r_j)h(r_k)\right]\\
        & = \frac{g(\tau) h(r_k)}{G} \frac{g(\tau) h(r_j)}{G} \frac{g(r_j) h(r_i)}{G}\\
        & = k(\tau, r_k) k(\tau, r_j) k(r_j, r_i).
    \end{align*}
\end{example}

Leveraging Definition ~\eqref{separable kernel} and Lemma ~\eqref{kernel product} we define an extended path \(\bar{\bX}:[0,T] \to \RR^{d'}\) such that the geometric Volterra signature \(\geovsig\) of the actual path \(\bX:[0,T] \to \RR^{d}\) can be realized as the classical signature of the extended path \(\bar{\bX}\). Here \(d'\) can be much larger than \(d\) depending on the context and the problem. In this construction, we consider \(d' = d +2\).

Let the augmented path \(\bar{\bX}:[0,T]\to\RR^{d+2}\)  be defined as 
\[\bar{\bX}_t:=\left(g(t),h(t),\bX_t^1,\ldots,\bX_t^d\right).\]
We associate with $\bar{\bX}$ the alphabet  \(\cA:=\left\{\mathbf{0},\bar{\mathbf{0}},\mathbf{1},\ldots,\mathbf{d}\right\},\) where \(\mathbf{0}\) and \(\bar{\mathbf{0}}\) correspond to the coordinates \(g\) and \(h\), respectively, and \(\mathbf{i}\) corresponds to the driving coordinate \(\bX^i\).

Using the standard coordinate basis, we identify \(T(\RR^{d+2})\) with the vector space of finite real linear combinations of words over \(\cA\). We denote the empty word by \(\emptyset\). Equipped with the shuffle product \(\shuffle\), this space is a commutative, associative, unital algebra with unit $\emptyset$. The shuffle product of two words is the sum of all interleavings that preserve the order of the letters within each word, counted with multiplicity.

\begin{thm}[Realization of \(\geovsig(\bX)\) as \(\sig(\Bar{\bX})\) for some \(\Bar{\bX}\)]
\label{classical realization of geovsig}
Under the preceding assumptions, let \(\vtree=(\ktree,\dtree)\in\cL_\nu\) have degree \(n\ge1\), and let \(m\) be the number of children of the dummy root in \(\ktree\). For every \(s\in[0,T)\), there exists a finite word polynomial \(\mathbf{P}_{\vtree}{(s)}\in T(\RR^{d+2})\), depending only on \(\vtree\), \(g(s)\), and \(h(s)\), such that
\[
    \innerprod{\geovsig(\bX)^\tau_{st},\vtree}
    =\frac{g(\tau)^m}{G^n}\innerprod{\sig(\bar{\bX})_{st},\mathbf{P}_{\vtree}^{(s)}},
    \qquad s<t\le\tau\le T.
\]
In particular, every geometric Volterra signature component is a finite linear combination of classical signature components of \(\bar{\bX}\).

Moreover, \(\mathbf{P}_{\vtree}{(s)}\) can be chosen to involve only words of lengths between \(n\) and \(3n-m\). Consequently, the classical signature of \(\bar{\bX}\) truncated at level \(3N-1\), together with the known kernel, determines every geometric Volterra signature component of degree at most \(N\).
\end{thm}

\begin{proof}
Using the definition of the geometric Volterra signature \(\geovsig\) of the path \(\bX\) and Lemma~\ref{kernel product} we have
    \begin{align} \label{gvsig reduction}
        \innerprod{\geovsig(\bX)^{\tau}_{st}, \vtree} & = \int_{s<r_1< \cdots < r_n < t} \left(\prod^n_{v=1} k (r_{\kparent{i_v}}, r_{i_v})\right) d\bX_{r_{i_1}}^{i_1} \cdots d\bX_{r_{i_n}}^{i_n} \nonumber\\
        & = \frac{1}{G^n}\int_{s<r_1< \cdots < r_n < t} \left(\prod^n_{v=1} g (r_{\kparent{i_v}}) h(r_{i_v})\right) d\bX_{r_{i_1}}^{i_1} \cdots d\bX_{r_{i_n}}^{i_n} \nonumber \\
        & = \frac{g(\tau)^m}{G^n}\int_{s<r_1< \cdots < r_n < t} \left(\prod^n_{v=1} g (r_{i_v})^{c_{i_v}} h(r_{i_v})\right) d\bX_{r_{i_1}}^{i_1} \cdots d\bX_{r_{i_n}}^{i_n}.
    \end{align}
Since
\begin{align}
   \prod^n_{v=1} \left(g (r_{i_v}) \right)^{c_{i_v}} & = \prod^n_{v=1} \left(g (r_{i_v}) - g(s) + g(s) \right)^{c_{i_v}} \nonumber \\
   & =  \prod^n_{v=1} \left(\innerprod{\sig(\bar{\bX})_{sr_{i_v}}, \mathbf{0}} + g(s) \right)^{c_{i_v}}\nonumber \\
   & =  \prod^n_{v=1} \left[\sum_{k_v = 0}^{c_{i_v}} \binom{c_{i_v}}{k_v} g(s)^{c_{i_v} - k_v}\left(\innerprod{\sig(\bar{\bX})_{sr_{i_v}}, \mathbf{0}}  \right)^{k_v} \right]\nonumber \\
   & = \sum_{k_1 = 0}^{c_{i_1}} \cdots \sum_{k_n = 0}^{c_{i_n}} \left( \prod_{v=1}^{n} \binom{c_{i_v}}{k_v}\right) g(s)^{\sum_{v=1}^{n} (c_{i_v} - k_v)} \prod_{v=1}^{n} \left(\innerprod{\sig(\bar{\bX})_{sr_{i_v}}, \mathbf{0}}  \right)^{k_v} \nonumber \\
   & = \sum_{k_1 = 0}^{c_{i_1}} \cdots \sum_{k_n = 0}^{c_{i_n}} \left( \prod_{v=1}^{n} \binom{c_{i_v}}{k_v}\right) g(s)^{\sum_{v=1}^{n} (c_{i_v} - k_v)} \prod_{v=1}^{n} \innerprod{\sig(\bar{\bX})_{sr_{i_v}}, \mathbf{0}^{\shuffle k_v}},
\end{align}
where
\[ \mathbf{0}^{\shuffle k_v} :=
\underbrace{\mathbf{0}\shuffle\cdots\shuffle\mathbf{0}
}_{k_v\text{ factors}} \quad\text{for }k_v\ge1, \qquad \mathbf{0}^{\shuffle 0}:=\mathbf{1}.
\]
Here $\mathbf{1}$ denotes the empty word, and the last
equality follows from the shuffle identity. Similarly,
\begin{align}
    \prod_{v=1}^{n} h(r_{i_v}) &= \prod_{v=1}^{n} \left(\innerprod{\sig(\bar{\bX})_{sr_{i_v}}, \bar{\mathbf{0}}}+h(s)\right) \nonumber\\
    &= \sum_{\ell_1=0}^{1}\cdots\sum_{\ell_n=0}^{1} h(s)^{n-\sum_{v=1}^{n}\ell_v} \prod_{v=1}^{n} \innerprod{\sig(\bar{\bX})_{sr_{i_v}}, \bar{\mathbf{0}}^{\shuffle\ell_v}}.
\end{align}
Let
\[
\cI_{\vtree}:=\left\{(\bk,\boldsymbol{\ell}): \bk=(k_1,\dots,k_n),\ \boldsymbol{\ell}=(\ell_1,\dots,\ell_n),\ k_v\in\{0,\ldots,c_{i_v}\},\ \ell_v\in\{0,1\}\ \text{for every }v=1,\ldots,n\right\}.
\]
For \((\mathbf{k},\boldsymbol{\ell})\in\mathcal I_{\vtree}\), define
\[
\alpha_{\mathbf{k},\boldsymbol{\ell}}(s):=\left(\prod_{v=1}^{n}\binom{c_{i_v}}{k_v}\right)g(s)^{\sum_{v=1}^{n}(c_{i_v}-k_v)}h(s)^{n-\sum_{v=1}^{n}\ell_v},
\qquad
L_{k_v,\ell_v}:=\mathbf{0}^{\shuffle k_v}\shuffle\bar{\mathbf{0}}^{\shuffle\ell_v},
\]
By the binomial theorem and the shuffle identity,
\begin{align} \label{kernel shuffle}
    \prod_{v=1}^{n}g(r_{i_v})^{c_{i_v}}h(r_{i_v}) &= \sum_{(\mathbf{k},\boldsymbol{\ell})\in\mathcal I_{\vtree}}\alpha_{\mathbf{k},\boldsymbol{\ell}}(s)\prod_{v=1}^{n}\innerprod{\sig(\bar{\bX})_{sr_{i_v}},L_{k_v,\ell_v}}.
\end{align}
Using these ~\eqref{kernel shuffle} in ~\eqref{gvsig reduction} gives
\begin{align} \label{geovsig using L}
        \innerprod{\geovsig(\bX)^{\tau}_{st}, \vtree}
        &= \frac{g(\tau)^m}{G^n}\sum_{(\mathbf{k},\boldsymbol{\ell})\in\mathcal I_{\vtree}} \alpha_{\mathbf{k},\boldsymbol{\ell}}(s)\int_{s<r_1<\cdots<r_n<t} \left(\prod_{v=1}^{n} \innerprod{\sig(\bar{\bX})_{sr_{i_v}},L_{k_v,\ell_v}}\right) d\bX_{r_{i_1}}^{i_1}\cdots d\bX_{r_{i_n}}^{i_n}.
    \end{align}
Fix \((\mathbf{k},\boldsymbol{\ell})\in\mathcal I_{\vtree}\) and write \(L_v:=L_{k_v,\ell_v}\). The innermost contribution in ~\eqref{geovsig using L} can simply be written as 
\[\int_s^{r_{i_2}} \innerprod{\sig(\bar{\bX})_{sr_{i_1}},L_1}\,d\bX_{r_{i_1}}^{i_1}=\innerprod{\sig(\bar{\bX})_{sr_{i_2}},L_1 i_1}.\]
The double integral contribution becomes
\begin{align*}
        & \int_{s<r_{i_1}<r_{i_2}<r_{i_3}} \innerprod{\sig(\bar{\bX})_{sr_{i_1}},L_1}\innerprod{\sig(\bar{\bX})_{sr_{i_2}},L_2}\,d\bX_{r_{i_1}}^{i_1}d\bX_{r_{i_2}}^{i_2}\\
        & \qquad = \int_s^{r_{i_3}} \innerprod{\sig(\bar{\bX})_{sr_{i_2}},L_1i_1}\innerprod{\sig(\bar{\bX})_{sr_{i_2}},L_2}\,d\bX_{r_{i_2}}^{i_2} \nonumber\\
        & \qquad = \innerprod{\sig(\bar{\bX})_{sr_{i_3}},((L_1i_1)\shuffle L_2)i_2}.
\end{align*}
Repeating this procedure gives
\[
        \int_{s<r_1<\cdots<r_n<t} \left(\prod_{v=1}^{n} \innerprod{\sig(\bar{\bX})_{sr_{i_v}},L_v}\right) d\bX_{r_{i_1}}^{i_1}\cdots d\bX_{r_{i_n}}^{i_n}
        =\innerprod{\sig(\bar{\bX})_{st},\bW_{\mathbf{k},\boldsymbol{\ell}}},
    \]
where is \(\bW_{\mathbf{k},\boldsymbol{\ell}}\) the finite word polynomial
    \[
        \bW_{\mathbf{k},\boldsymbol{\ell}}:=\bigl(\cdots\bigl(((L_1i_1)\shuffle L_2)i_2\shuffle L_3\bigr)i_3\cdots\shuffle L_n\bigr)i_n, \quad \text{ and } \quad \bW_{\mathbf{k},\boldsymbol{\ell}}:=L_1i_1, \text{ for } n=1.
    \]
Substituting this identity in ~\eqref{geovsig using L} yields  
 \[
    \innerprod{\geovsig(\bX)^\tau_{st},\vtree}=\frac{g(\tau)^m}{G^n}\sum_{(\mathbf{k},\boldsymbol{\ell})\in\mathcal I_{\vtree}} \alpha_{\mathbf{k},\boldsymbol{\ell}}(s)\innerprod{\sig(\bar{\bX})_{st},W_{\mathbf{k},\boldsymbol{\ell}}}.
    \]
    Each \(\bW_{\mathbf{k},\boldsymbol{\ell}}\) is a finite linear combination of words. If we define
\[ \mathbf{P}_{\vtree}{(s)}:=\sum_{(\mathbf{k},\boldsymbol{\ell})\in\cI_{\vtree}} \alpha_{\mathbf{k},\boldsymbol{\ell}}(s)\mathbf{W}_{\mathbf{k},\boldsymbol{\ell}}.
\]
Then the preceding calculation proves the first part of the claim. Now, each word occurring in \(\bW_{\mathbf{k},\boldsymbol{\ell}}\) has length
\[
    n+\sum_{v=1}^{n}(k_v+\ell_v)\le n+\sum_{v=1}^{n}c_v+n=3n-m.
\]
Since $m\ge1$ for every nonempty kernel tree, $3n-m\le3N-1$ whenever $n\le N$.

\end{proof}

\bibliographystyle{plain}
\bibliography{myrefs}

@book{Friz_Victoir_2010, place={Cambridge}, series={Cambridge Studies in Advanced Mathematics}, title={Multidimensional Stochastic Processes as Rough Paths: Theory and Applications}, publisher={Cambridge University Press}, author={Friz, Peter K. and Victoir, Nicolas B.}, year={2010}, collection={Cambridge Studies in Advanced Mathematics}}

@article{cuchiero2023signature,
  title={Signature-based models: theory and calibration},
  author={Cuchiero, Christa and Gazzani, Guido and Svaluto-Ferro, Sara},
  journal={SIAM journal on financial mathematics},
  volume={14},
  number={3},
  pages={910--957},
  year={2023},
  publisher={SIAM}
}

@article{lyons1998differential,
  title={Differential equations driven by rough signals},
  author={Lyons, Terry J},
  journal={Revista Matem{\'a}tica Iberoamericana},
  volume={14},
  number={2},
  pages={215--310},
  year={1998}
}

@incollection{connes1999hopf,
  title={Hopf algebras, renormalization and noncommutative geometry},
  author={Connes, Alain and Kreimer, Dirk},
  booktitle={Quantum field theory: perspective and prospective},
  pages={59--109},
  year={1999},
  publisher={Springer}
}

@article{manchon2008hopf,
  title={Hopf algebras in renormalisation},
  author={Manchon, Dominique},
  journal={Handbook of algebra},
  volume={5},
  pages={365--427},
  year={2008},
  publisher={Elsevier}
}

@book{abe2004hopf,
  title={Hopf algebras},
  author={Abe, Eiichi},
  volume={74},
  year={2004},
  publisher={Cambridge university press}
}

@article{brouder2004trees,
  title={Trees, renormalization and differential equations},
  author={Brouder, Ch},
  journal={BIT Numerical Mathematics},
  volume={44},
  number={3},
  pages={425--438},
  year={2004},
  publisher={Springer}
}

@article{vl_ComteRenault1996,
  author = {Comte, F. and Renault, E.},
  title = {{Long memory continuous time models}},
  year = {1996},
  journal = {Journal of Econometrics},
  volume = {73},
  number = {1},
  pages = {101--149},
}

@article{vl_Gatheral2018,
  author = {Gatheral, Jim and Jaisson, Thibault and Rosenbaum, Mathieu},
  title = {{Volatility is rough}},
  year = {2018},
  journal = {Quantitative Finance},
  volume = {18},
  number = {6},
  pages = {933--949},
}

@article{vl_GuyonLekeufack2023,
  author = {Guyon, Julien and Lekeufack, Jordan},
  title = {{Volatility is (mostly) path-dependent}},
  year = {2023},
  journal = {Quantitative Finance},
  volume = {23},
  number = {9},
  pages = {1221--1258},
}

@article{vl_BoydChua1985,
  author = {Boyd, S. and Chua, L.},
  title = {{Fading memory and the problem of approximating nonlinear operators with Volterra series}},
  year = {1985},
  journal = {IEEE Transactions on Circuits and Systems},
  volume = {32},
  number = {11},
  pages = {1150--1161},
}

@article{vl_Leland1994,
  author = {Leland, W. E. and Taqqu, M. S. and Willinger, W. and Wilson, D. V.},
  title = {{On the self-similar nature of Ethernet traffic (extended version)}},
  year = {1994},
  journal = {IEEE/ACM Transactions on Networking},
  volume = {2},
  number = {1},
  pages = {1--15},
}

@article{vl_BurovBarkai2008,
  author = {Burov, S. and Barkai, E.},
  title = {{Fractional Langevin equation: Overdamped, underdamped, and critical behaviors}},
  year = {2008},
  journal = {Physical Review E},
  volume = {78},
  number = {3},
  pages = {031112},
}

@incollection{vl_CarmonaCoulon2014,
  author = {Carmona, Ren{\'{e}} and Coulon, Michael},
  title = {{A survey of commodity markets and structural models for electricity prices}},
  year = {2014},
  booktitle = {Quantitative Energy Finance},
  pages = {41--83},
  publisher = {Springer},
}

@article{alqahtani2023solution,
  title={Solution of ill-posed problems with Chebfun},
  author={Alqahtani, Abdulaziz and Mach, Thomas and Reichel, Lothar},
  journal={Numerical Algorithms},
  volume={92},
  number={4},
  pages={2341--2364},
  year={2023},
  publisher={Springer}
}

@incollection{vl_DiethelmFreed1999,
  author = {Diethelm, Kai and Freed, Alan D.},
  title = {{On the solution of nonlinear fractional-order differential equations used in the modeling of viscoplasticity}},
  year = {1999},
  booktitle = {Scientific Computing in Chemical Engineering II},
  pages = {217--224},
  publisher = {Springer},
  doi = {10.1007/978-3-642-60185-9_24},
}

@article{vl_Cromer1985,
  author = {Cromer, Thomas L.},
  title = {{Asymptotically periodic solutions to Volterra integral equations in epidemic models}},
  year = {1985},
  journal = {Journal of Mathematical Analysis and Applications},
  volume = {110},
  number = {2},
  pages = {483--494},
  doi = {10.1016/0022-247X(85)90310-5},
}

@book{vl_GerstnerKistler2002,
  author = {Gerstner, Wulfram and Kistler, Werner M.},
  title = {{Spiking Neuron Models: Single Neurons, Populations, Plasticity}},
  year = {2002},
  publisher = {Cambridge University Press},
}

@article{vl_Kou2008,
  author = {Kou, S. C.},
  title = {{Stochastic modeling in nanoscale biophysics: Subdiffusion within proteins}},
  year = {2008},
  journal = {The Annals of Applied Statistics},
  volume = {2},
  number = {2},
  pages = {501--535},
}

@article{vl_Elman1990,
  author = {Elman, Jeffrey L.},
  title = {{Finding structure in time}},
  year = {1990},
  journal = {Cognitive Science},
  volume = {14},
  number = {2},
  pages = {179--211},
}

@incollection{vl_Jordan1997,
  author = {Jordan, Michael I.},
  title = {{Serial order: A parallel distributed processing approach}},
  year = {1997},
  booktitle = {Neural-Network Models of Cognition: Biobehavioral Foundations},
  volume = {121},
  pages = {471--495},
  publisher = {Elsevier Science Press},
  series = {Advances in Psychology},
  doi = {10.1016/S0166-4115(97)80111-2},
}

@article{vl_Hochreiter1997,
  author = {Hochreiter, Sepp and Schmidhuber, J{\"{u}}rgen},
  title = {{Long short-term memory}},
  year = {1997},
  journal = {Neural Computation},
  volume = {9},
  number = {8},
  pages = {1735--1780},
}

@inproceedings{vl_Vaswani2017,
  author = {Vaswani, Ashish and Shazeer, Noam and Parmar, Niki and Uszkoreit, Jakob and Jones, Llion and Gomez, Aidan N. and Kaiser, {\L}ukasz and Polosukhin, Illia},
  title = {{Attention is all you need}},
  year = {2017},
  booktitle = {Advances in Neural Information Processing Systems},
}

@misc{vl_Wang2020,
  author = {Wang, Sinong and Li, Belinda Z. and Khabsa, Madian and Fang, Han and Ma, Hao},
  title = {{Linformer: Self-attention with linear complexity}},
  year = {2020},
  eprint = {2006.04768},
  archivePrefix = {arXiv},
  url = {https://arxiv.org/abs/2006.04768},
  howpublished = {arXiv preprint arXiv:2006.04768},
}

@inproceedings{vl_Gu2022,
  author = {Gu, Albert and Goel, Karan and R{\'{e}}, Christopher},
  title = {{Efficiently modeling long sequences with structured state spaces}},
  year = {2022},
  booktitle = {International Conference on Learning Representations},
}

@misc{vl_GuDao2023,
  author = {Gu, Albert and Dao, Tri},
  title = {{Mamba: Linear-time sequence modeling with selective state spaces}},
  year = {2023},
  eprint = {2312.00752},
  archivePrefix = {arXiv},
  url = {https://arxiv.org/abs/2312.00752},
  howpublished = {arXiv preprint arXiv:2312.00752},
}

@inproceedings{vl_DaoGu2024,
  author = {Dao, Tri and Gu, Albert},
  title = {{Transformers are SSMs: Generalized models and efficient algorithms through structured state space duality}},
  year = {2024},
  booktitle = {Proceedings of the 41st International Conference on Machine Learning},
  pages = {10041--10071},
  publisher = {PMLR},
}

@inproceedings{vl_Chen2018,
  author = {Chen, Ricky T. Q. and Rubanova, Yulia and Bettencourt, Jesse and Duvenaud, David},
  title = {{Neural ordinary differential equations}},
  year = {2018},
  booktitle = {Advances in Neural Information Processing Systems},
  volume = {31},
  pages = {6572--6583},
  publisher = {Curran Associates},
}

@inproceedings{vl_Kidger2020,
  author = {Kidger, Patrick and Morrill, James and Foster, James and Lyons, Terry},
  title = {{Neural controlled differential equations for irregular time series}},
  year = {2020},
  booktitle = {Advances in Neural Information Processing Systems},
  volume = {33},
  url = {https://papers.nips.cc/paper/2020/hash/4a5876b450b45371f6cfe5047ac8cd45-Abstract.html},
}

@article{vl_Zappala2024,
  author = {Zappala, Emanuele and {de Oliveira Fonseca}, Antonio Henrique and {Ortega Caro}, Josue and Moberly, Andrew Henry and Higley, Michael James and Cardin, Jessica and {van Dijk}, David},
  title = {{Learning integral operators via neural integral equations}},
  year = {2024},
  journal = {Nature Machine Intelligence},
  volume = {6},
  number = {9},
  pages = {1046--1062},
}

@article{vl_Bergerhausen2025,
  author = {Bergerhausen, Martin and Pr{\"{o}}mel, David J. and Scheffels, David},
  title = {{Neural stochastic Volterra equations: Learning path-dependent dynamics}},
  year = {2025},
  journal = {Journal of Machine Learning},
  volume = {4},
  number = {4},
  pages = {264--289},
}

@article{vl_Lipton2018,
  author = {Lipton, Zachary C.},
  title = {{The mythos of model interpretability}},
  year = {2018},
  journal = {Communications of the ACM},
  volume = {61},
  number = {10},
  pages = {36--43},
}

@article{vl_Rudin2019,
  author = {Rudin, Cynthia},
  title = {{Stop explaining black box machine learning models for high stakes decisions and use interpretable models instead}},
  year = {2019},
  journal = {Nature Machine Intelligence},
  volume = {1},
  pages = {206--215},
}

@article{vl_Bengio1994,
  author = {Bengio, Yoshua and Simard, Patrice and Frasconi, Paolo},
  title = {{Learning long-term dependencies with gradient descent is difficult}},
  year = {1994},
  journal = {IEEE Transactions on Neural Networks},
  volume = {5},
  number = {2},
  pages = {157--166},
}

@inproceedings{vl_Pascanu2013,
  author = {Pascanu, Razvan and Mikolov, Tomas and Bengio, Yoshua},
  title = {{On the difficulty of training recurrent neural networks}},
  year = {2013},
  booktitle = {Proceedings of the 30th International Conference on Machine Learning},
  pages = {1310--1318},
}

@article{vl_Raissi2019,
  author = {Raissi, Maziar and Perdikaris, Paris and Karniadakis, George E.},
  title = {{Physics-informed neural networks: A deep learning framework for solving forward and inverse problems involving nonlinear partial differential equations}},
  year = {2019},
  journal = {Journal of Computational Physics},
  volume = {378},
  pages = {686--707},
}

@article{vl_Karniadakis2021,
  author = {Karniadakis, George Em and Kevrekidis, Ioannis G. and Lu, Lu and Perdikaris, Paris and Wang, Sifan and Yang, Liu},
  title = {{Physics-informed machine learning}},
  year = {2021},
  journal = {Nature Reviews Physics},
  volume = {3},
  number = {6},
  pages = {422--440},
}

@book{vl_Corduneanu2008,
  author = {Corduneanu, Constantin},
  title = {Integral Equations and Applications},
  year = {2008},
  publisher = {Cambridge University Press},
  note = {Reprint of the 1991 edition},
}

@book{vl_Gripenberg1990,
  author = {Gripenberg, Gustaf and Londen, Stig-Olof and Staffans, Olof},
  title = {Volterra Integral and Functional Equations},
  year = {1990},
  publisher = {Cambridge University Press},
  series = {Encyclopedia of Mathematics and its Applications},
  volume = {34},
}

@article{vl_PardouxProtter1990,
  author = {Pardoux, Etienne and Protter, Philip},
  title = {Stochastic {Volterra} equations with anticipating coefficients},
  year = {1990},
  journal = {The Annals of Probability},
  volume = {18},
  number = {4},
  pages = {1635--1655},
}

@article{harang2024volterra,
  author = {Harang, Fabian and Tindel, Samy and Wang, Xiaohua},
  title = {{Volterra equations driven by rough signals 3:
    Probabilistic construction of the Volterra rough path
    for fractional Brownian motions}},
  journal = {Journal of Theoretical Probability},
  year = {2024},
  volume = {37},
  pages = {307--351},
  doi = {10.1007/s10959-023-01251-y},
}

@incollection{vl_OksendalZhang1993,
  author = {{\O}ksendal, Bernt and Zhang, Tu-Sheng},
  title = {The stochastic {Volterra} equation},
  year = {1993},
  booktitle = {The Barcelona Seminar on Stochastic Analysis},
  editor = {Nualart, David and Sanz-Sol{\'e}, Marta},
  series = {Progress in Probability},
  volume = {32},
  pages = {168--202},
  publisher = {Birkh{\"a}user},
}

@article{ali2025branched,
  title={Branched Signature Model},
  author={Ali, Munawar and Feng, Qi},
  journal={arXiv preprint arXiv:2511.00018},
  year={2025}
}

@article{harang2021volterra,
  title={Volterra equations driven by rough signals},
  author={Harang, Fabian A and Tindel, Samy},
  journal={Stochastic Processes and their Applications},
  volume={142},
  pages={34--78},
  year={2021},
  publisher={Elsevier}
}

@article{gubinelli2010ramification,
  title={Ramification of rough paths},
  author={Gubinelli, Massimiliano},
  journal={Journal of Differential Equations},
  volume={248},
  number={4},
  pages={693--721},
  year={2010},
  publisher={Elsevier}
}

@article{gubinelli2004controlling,
  title={Controlling rough paths},
  author={Gubinelli, Massimiliano},
  journal={Journal of Functional Analysis},
  volume={216},
  number={1},
  pages={86--140},
  year={2004},
  publisher={Elsevier}
}

@article{harang2023volterra,
  title={Volterra equations driven by rough signals 2: Higher-order expansions},
  author={Harang, Fabian A and Tindel, Samy and Wang, Xiaohua},
  journal={Stochastics and Dynamics},
  volume={23},
  number={01},
  pages={2350002},
  year={2023},
  publisher={World Scientific}
}

@article{chen1977iterated,
  title={Iterated path integrals},
  author={Chen, Kuo-Tsai},
  journal={Bulletin of the American Mathematical Society},
  volume={83},
  number={5},
  pages={831--879},
  year={1977}
}

@article{feng2023cubature,
  title={Cubature Method for Stochastic Volterra Integral Equations},
  author={Feng, Qi and Zhang, Jianfeng},
  journal={SIAM Journal on Financial Mathematics},
  volume={14},
  number={4},
  pages={959--1003},
  year={2023},
  publisher={SIAM}
}

@incollection{chevyrev2025primer,
  title={A primer on the signature method in machine learning},
  author={Chevyrev, Ilya and Kormilitzin, Andrey},
  booktitle={Signature Methods in Finance: An Introduction with Computational Applications},
  pages={3--64},
  year={2025},
  publisher={Springer Nature Switzerland, Cham}
}

@misc{castrequini2014youngsystems,
      title={On Young Systems}, 
      author={R. A. Castrequini and P. J. Catuogno},
      year={2014},
      eprint={1412.1970},
      archivePrefix={arXiv},
      primaryClass={math.AP},
      url={https://arxiv.org/abs/1412.1970}, 
}

@article{bruned2023ramification,
  title={Ramification of Volterra-type rough paths},
  author={Bruned, Yvain and Katsetsiadis, Foivos},
  journal={Electronic Journal of Probability},
  volume={28},
  pages={1--25},
  year={2023},
  publisher={The Institute of Mathematical Statistics and the Bernoulli Society}
}

@article{hager2026volterra,
  title={The Volterra signature},
  author={Hager, Paul P and Harang, Fabian N and Pelizzari, Luca and Tindel, Samy},
  journal={arXiv preprint arXiv:2603.04525},
  year={2026}
}

@article{cass2022combinatorial,
  title={A combinatorial approach to geometric rough paths and their controlled paths},
  author={Cass, Thomas and Driver, Bruce K and Litterer, Christian and Rossi Ferrucci, Emilio},
  journal={Journal of the London Mathematical Society},
  volume={106},
  number={2},
  pages={936--981},
  year={2022},
  publisher={Wiley Online Library}
}
\end{document}